\documentclass{amsart}

\usepackage{amsmath,amssymb,amsbsy,amsfonts, latexsym, amsopn,amstext, amsxtra,euscript,amscd}   

\usepackage[colorlinks=true, urlcolor=blue, citecolor=blue,  linkcolor=blue, pdfborder={0 0 0}]{hyperref}

\usepackage[msc-links]{amsrefs}

\usepackage{tikz-cd}
   
\input xy \xyoption{all}

\newtheorem{thm}{Theorem}
\newtheorem{prop}{Proposition}
\newtheorem{lem}{Lemma}
\newtheorem{rem}{Remark}
\newtheorem{cor}{Corollary}
\newtheorem{defn}{Definition}
\newtheorem{exa}{Example}

\usepackage[capitalise]{cleveref}

\crefname{thm}{Thm.}{}
\crefname{prop}{Prop.}{}
\crefname{lem}{Lem.}{}
\crefname{cor}{Cor.}{}

\DeclareMathOperator\chara{char }
\DeclareMathOperator\J{J}
\DeclareMathOperator\bAut{\overline{Aut}\, }
\DeclareMathOperator\Aut{Aut \, }

\DeclareMathOperator\Pic{Pic }
\DeclareMathOperator\Jac{Jac }

\DeclareMathOperator\Kum{Kum }

\DeclareMathOperator\Hom{Hom }

\DeclareMathOperator\Img{Img }
\DeclareMathOperator\Spec{Spec}

\DeclareMathOperator\KT{K3 }
\DeclareMathOperator\SI{SI }

\def\Q{\mathbb Q}
\def\bP{\mathbb P}
\def\P{\mathbb P}

\def\Z{\mathbb Z}

\def\F{\mathbb F}

\newcommand\A{\mathcal A}
\newcommand\B{\mathcal B}
\newcommand\X{\mathcal C}
\newcommand\DD{\mathcal D}
\newcommand\M{\mathcal M}
\newcommand\Y{\mathcal Y}
\newcommand\E{\mathcal E}

\def\H{\mathcal H}
\def\L{\mathcal L}
\def\S{\mathcal S}
\def\O{\mathcal O}

\def\u{\mathfrak u}
\def\v{\mathfrak v}
\def\p{\mathfrak p}

\newcommand\e{\varepsilon}

\newcommand\bG{\overline{G}}

\def\<{\langle}
\def\>{\rangle}

\def\iso{{\, \cong\, }}

\title{Isogenous components of Jacobian surfaces}
\thanks{This is a corrected version of the article published in Eur. J. Math. \textbf{6} (2020), no.~4, 1276--1302, \url{https://doi.org/10.1007/s40879-019-00375-y}. The published version contains errors in genus computations computed via Maple, which resulted in errors  in  \cref{prop-2},   \cref{prop-4}, and \cref{thm-1}. Several components of the isogeny loci were reported there with geometric genus at least two, while they have genus zero.  As a consequence, some finiteness statements of the published version are false. }

\author{L. Beshaj}
\address{Depart.  of Mathematical Sciences\\
United States Military Academy at West Point
}
 \email{Lubjana.Beshaj@usma.edu}

\author{A. Elezi}
\address{Depart. of Mathematics and Statistics \\
American University, \\
Washington DC, 20016.
}
\email{aelezi@american.edu}

\author{T. Shaska}
\address{Depart. of Mathematics and Statistics\\
Oakland University
Rochester, MI, 48309.
}
\email{shaska@oakland.edu}

\begin{document}

\begin{abstract}
Let $\X$ be a genus 2 curve defined over a  field $K$, $\chara K = p \geq 0$, and $\Jac (\X, \iota)$ its Jacobian, where $\iota$ is the principal polarization of $\Jac (\X)$ attached to $\X$. Assume that $\Jac (\X)$ is geometrically $(n,n)$-reducible with  $E_1$ and $E_2$  its elliptic components.  When $K$ is a number field, we count such curves up to geometric isomorphism, among those admitting a model over $K$. For $n=2$ and $\Aut(\X)\iso V_4$ we prove that, for each $N=2,3,5,7$, there are infinitely many curves $\X$ for which $E_1$ and $E_2$ are $N$-isogenous. For $n=2$ and $\Aut(\X)\iso D_4$ the elliptic quotients from the two classes of nonhyperelliptic involutions are always $2$-isogenous, and for each $N=3,5,7$ only finitely many curves have them $N$-isogenous. For $n=3$ there are infinitely many curves with $N$-isogenous elliptic components for each $N=2,3,5$, while for $N=7$ there are only finitely many such curves with a model over $\Q$.  Furthermore, we determine the Kummer  and the Shioda-Inose surfaces for the above $\Jac \X$ and show how such  results  in positive characteristic $p>2$ suggest nice applications in cryptography. 
\end{abstract}

\maketitle

\setcounter{tocdepth}{1}


\section{Introduction}
An Abelian variety $\A$, defined over a field $k$,   is   simple if it has no proper non-zero Abelian subvariety over $k$. $\A$ is called reducible (or decomposable) if it is isogenous to a direct product of Abelian varieties.  We call $\A$ geometrically simple (or absolutely simple) if it is simple over the algebraic closure of $k$,  analogously   $\A$ is geometrically reducible when it is reducible over the algebraic closure of $k$.  
%
A 2-dimensional Jacobian variety is geometrically reducible if and only if it is $(n, n)$-decomposable for some $n>1$.  Reducible Jacobian varieties have been studied   extensively since the XIX-century, most notably by Fricke, Clebsch, and Bolza.    In the late XX-century they became the focus of many mathematicians through the work of Frey \cite{frey-95, frey-kani},    Shaska and  V\"olklein \cite{deg2, deg3, deg5}, Kumar \cite{kumar} and many others.  
If $\A/k$ is a 2-dimensional reducible Jacobian  variety defined over a field $k$,   then there is a degree $n^2$  isogeny to a product $\A \iso E_1 \times E_2$, where $E_i$, $i=1, 2$ are 1-dimensional. The main focus of this paper is to investigate when $E_1$ and $E_2$  are  isogenous to each other and how often this occurs for a fixed $n$.  

The question has received attention lately for different reasons.  In \cite{kumar-kuwata} the authors were able to determine the Mordell-Weil rank of elliptic fibrations $F^{(i)}$, for $i=1, \dots , 6$;    when $E_1$ and $E_2$ were isogenous, not isomorphic, and with complex multiplication, and show that in this case both $F^5$ and $F^6$ have rank 18.      
%
In recent   developments in supersingular isogeny based cryptography (SIDH) Costello \cite{costello} focuses on the $(2, 2)$ reducible Jacobians, where the addition  is done via the Kummer surface.   More importantly, it seems that the most interesting case is exactly the case when $E_1$ is isogenous to $E_2$. In this case, since the decomposition of the Abelian varieties is determined up to isogeny, the 2-dimensional Jacobian is isogenous to $E^2$.  There are several interesting questions that arise when we consider such Jacobians over the finite field $\F_p$.  

The focus of this paper is to investigate when the two elliptic components of the $(n, n)$ reducible 2-dimensional Jacobians are isogenous to each other.  
The space of genus 2 curves with $(n, n)$ reducible   Jacobians,  for  $n=2$ or $n$ is odd is a 2-dimensional irreducible locus $\L_n$ in the moduli space of curves $\M_2$.  For $n=2$ this   is the well known locus of curves with extra involutions \cite{deg2, m-sh, m-sh-3}, for $n$  odd such spaces were computed for the first time in \cite{thesis}, \cite{deg3}, \cite{deg5}. 
If $E_1$ and $E_2$ are $N$-isogenous then their $j$-invariants $j_1$ and $j_2$ satisfy the equation of the modular curve $X_0(N)$, say $\mathcal S_N := \phi_N (j_1, j_2)=0$. Such a curve can be embedded in $\M_2$.  So we want to study the intersection between $\L_n$ and $\mathcal S_N$ for given $n$ and $N$.  More precisely, for any number field $K$ we want to determine the number of $K$-rational points of this intersection. 

Our approach is computational.  We will focus on the cases when $n=2, 3$ and $N=2, 3, 5, 7$.   We prove that for $n=2$ and $N=2,3,5,7$ there are infinitely many curves $\X$ defined over $K$ with $\Aut (\X)\iso V_4$ such that $E_1$ and $E_2$ are $N$-isogenous; if $\Aut (\X)\iso D_4$ (resp. $D_6$) there is a 1-dimensional family such that $E_1$ and $E_2$ are 2-isogenous (resp. 3-isogenous), and for $\Aut (\X)\iso D_4$ and $N=3,5,7$ there are only finitely many such curves. For $n=3$ there are infinitely many curves $\X$ defined over $K$ such that $E_1, E_2$ are $N$-isogenous for $N=2,3,5$, while for $N=7$ there are only finitely many such curves defined over $\Q$.   The finiteness statements follow from explicit equations: in the $D_4$ case the isogeny condition is a nonzero polynomial in one parameter, and for $(n,N)=(3,7)$ we use a rank computation and a real obstruction.

Our paper is organized as follows. In \cref{sect-2} we give a brief account of the basic definitions of Abelian varieties and their isogenies. 
In \cref{sect-4} we first prove that for $n=2$ there are infinitely many genus $2$ curves $\X$ defined over a number field  $K$ with $\Aut (\X)\iso V_4$ whose elliptic components are $N$-isogenous for $N=2,3,5,7$.  Also,  for $n=2$ and $N = 3,5, 7$, there are only finitely many such $\X$ (up to isomorphism) with $\Aut(\X) = D_4$.  That $\X$ is defined over $K$ follows from the important fact that the invariants $u$ and $v$ are in the field of moduli of the curve $\X$ and that for every curve in $\L_2$, the field of moduli is a field of definition; see \cite{m-sh}.  This is not necessarily true for curves in $\L_n$, when $n>2$. However, a proof of the above result is still possible using the computational approach by using invariants $r_1$, $r_2$ of two cubics in \cite{deg3}.  These invariants are denoted by $\chi$ and $\psi$ here. 

Then we study the $n=3$ case. The equation of $\L_3$ was computed by the third author in his thesis \cite{thesis} and summarized  in \cite{deg3}. 
A birational parametrization of $\L_3$ exists in terms of the invariants $\chi$ and $\psi$.   We are able to compute the $j$-invariants of $E_1$ and $E_2$ in terms of $\chi$ and $\psi$ and find the conditions that $\chi$ and $\psi$ must satisfy. 
 Since  ordered  pairs $(\chi, \psi)$ are in one-to-one correspondence with genus two curves with $(3, 3)$-split Jacobians, then we  determine pairs $(\chi, \psi)$ such that the corresponding $j$-invariants $j_1$ and $j_2$  satisfy the equation of the modular curve $X_0 (N)$.  This case is different from $n=2$ in that a rational ordered pair $(\chi, \psi)$ does not necessarily correspond to a genus two curve defined over $K$.  However, a genus two curve defined over $K$ gives rise to rational invariants $\chi, \psi \in K$.  Hence, it is enough to count the rational ordered pairs $(\chi, \psi)$ that satisfy the equation $\phi_N (j_1, j_2   )=0$   of the modular curve $X_0 (N)$.        We prove that for $N=5$ there are infinitely many genus two curves $\X$ defined over $K$ such that they have $(3, 3)$-split Jacobian and $E_1$ and $E_2$ are $5$-isogenous: the curve $\phi_5(j_1, j_2)=0$ has a rational component.   
For $N=2,3$ there are also infinitely many such curves, while for $N=7$ there are only finitely many such curves admitting models over $\Q$.

 In the last section we consider the Kummer and Shioda-Inose surfaces of $(n, n)$-reducible Jacobians.  We classify such surfaces  when $n=2, 3$ (cf. \cref{prop-10}  and \cref{cor-2}). 

All computations in this paper were carried out with the genus $2$ package \cite{genus2}; the computations specific to this paper accompany its arXiv version as ancillary files.

\section{Preliminaries}\label{sect-2}  

An Abelian variety defined over $k$ is an absolutely irreducible projective variety defined over $k$ which is  a group scheme.  We will denote an Abelian variety  defined over a field $k$ by $\A_k$ or simply $\A$ when there is no confusion.   A morphism from the  Abelian variety $\A_1$ to the Abelian variety $\A_2$ is a homomorphism if and only if it maps the identity element of $\A_1$ to the identity element of $\A_2$.  

An abelian variety over a field $k$ is called \textit{simple} if it has no proper non-zero Abelian subvariety over $k$,  it is called \textit{absolutely simple} (or \textit{geometrically simple}) if it   is simple over the algebraic closure of $k$.    An Abelian variety of dimension 1 is called an \textit{elliptic curve}.

%
A homomorphism $f : \A \to  \B$ is called an \textit{isogeny} if $\Img f = \B$  and $\ker f$ is a finite group scheme. If an isogeny $\A \to \B$ exists we say that $\A$ and $\B$ are isogenous.  This relation is symmetric.  The degree of an isogeny $f : \A \to \B$ is the degree of the function field extension
$ \deg f  := [k(\A) : f^\star k(\B)]$.   It is equal to the order of the group scheme $\ker (f)$, which is, by definition, the scheme theoretical inverse image $f^{-1}(\{0_\A\})$.

The group of $\bar{k}$-rational points has order $\#(\ker f)(\bar{k}) = [k(A) : f^\star k(B)]^{sep}$, where $[k(A) : f^\star k(B)]^{sep}$ is the degree of the maximally separable extension in $k(\A)/ f^\star k(\B)$.  We say that   $f$ is a \textit{separable isogeny} if and only if  $\# \ker f(\bar{k}) = \deg f$.

For any Abelian variety $\A/k$ there is a one to one correspondence between the finite subgroup schemes $H \leq \A$ and   isogenies $f : \A \to \B$, where $\B$ is determined up to isomorphism.   Moreover, $H = \ker f$ and $\B = \A/H$. 
$f$ is separable if and only if $H$ is   \'{e}tale, and then $\deg f =\# H (\bar{k})$.
%
The following is often called the fundamental theorem of Abelian varieties.   Let $\A$ be an Abelian variety.  Then $\A$ is isogenous to 
\[   \A_1^{n_1} \times \A_2^{n_2} \times \dots \times \A_r^{n_r}, \]
where (up to permutation of the factors) $\A_i$ , for $i=1, \dots , r$ are simple, non-isogenous, Abelian varieties.  Moreover,  up to permutations,  the factors  $\A_i^{n_i}$ are uniquely determined  up to isogenies. 

When $k=\bar{k}$, then  let $f$ be a nonzero  isogeny   of $\A$. Its kernel $\ker f$ is a subgroup scheme of $\A$. It contains $0_\A$ and so its connected component, which is, by definition, an Abelian variety. 

%
Let $\X$ be a curve of genus 2 defined over a perfect field $k$ such that $\chara  k \neq 2$ and $\J = \Jac (\X)$ its Jacobian.  Fix a prime $\ell \geq 3$ and let $S$ be a maximal $\ell$-Weil isotropic subgroup of $\J [\ell]$, then we have $S\iso (\Z/\ell\Z)^2$.  Let $\J^\prime := \J/S$ be the quotient variety and $\Y$ a genus 2 curve such that $\Jac (\Y)=\J^\prime$. Hence, the classical isogeny problem becomes to compute $\Y$ when given $\X$ and $S$. 
 
If $\ell =2$ this problem is done with the Richelot construction.  Over finite fields this is done by Lubicz and Robert in \cite{L-R-2} using theta-functions. In general, if $\phi : \J (\X) \to \J (\Y)$ is the isogeny and $\Theta_\X$, $\Theta_Y$ the corresponding theta divisors, then $\phi (\Theta_\X)$ is in $| \ell \Theta_\Y|$. Thus, the image of $\phi (\Theta_\X)$ in the Kummer surface $\mathcal K_\Y = \J (\Y)/\<\pm 1\>$ is a degree $2\ell$ genus zero curve in $\P^3$ of arithmetic genus $\frac 1 2 \, (\ell^2-1)$. This curve can be computed without knowing $\phi$; see \cite{D-L} or \cite{frey-shaska} for details.   

For $\X$ given by $y^2=f(x)$, we have the divisor at infinity
\[ D_\infty := (1 : \sqrt{f(x)} : 0) + (1 : - \sqrt{f(x)} : 0)  \] 
The Weierstrass points of $\X$ are the projective roots of $f(x)$, namely $w_i:=(x_i, z_i)$, for $i=1, \dots , 6$  and  the Weierstrass divisor $W_\X$ is 
$ W_\X := \sum_{i=1}^6  (x_i, 0, z_i)$. 
A canonical divisor on $\X$ is $ \mathcal K_\X= W_\X - 2 D_\infty$. 
Let $D \in \Jac \X$ be a divisor expressed as $D = P+Q- D_\infty$.  The effective divisor $P+Q$ is determined by an ideal of the form $(a(x), y-b(x))$, where $b(x)$ is a cubic and $a(x)$ a monic polynomial of degree $d\leq 2$.

%

We can define the $\ell$-tuple embedding  $\rho_{2\ell} : \P^2 \to \P^{2\ell}$ by 
\[ (x, y, z) \to ( z^{2\ell}, \dots , x^i z^{2\ell -i}, x^{2\ell}    )\]
and denote the image of this map by $\mathcal R_{2\ell}$.  It is a rational normal curve of degree $2\ell$ in $\P^{2\ell}$.  Hence, any $2\ell +1$ distinct points on $\mathcal R_{2\ell}$ are linearly independent. Therefore, the images under $\rho_{2\ell}$ of the Weierstrass points of $\X$ are linearly independent for $\ell \geq 3$. Thus, the subspace $W := \< \rho_{2\ell} (W_\X) \> \subset \P^{2\ell}$    is $5$-dimensional.    For any pair of points $P, Q$ in  $\X$, the secant line $\mathcal L_{P, Q}$ is defined to be the line in $\P^{2\ell}$ intersecting $\mathcal R_{2\ell}$ in $\rho_{2\ell} (P) + \rho_{2\ell} (Q)$. In other words, 
\[
\mathcal L_{P, Q} = \left\{
\begin{split}
& \< \rho_{2\ell} (P), \rho_{2\ell} (Q) \> \; \; \text{ if } P \not\in  \{ Q, \tau (Q)\}  \\
& T_{\rho_{2\ell} (P)} (\mathcal R_{2\ell}) \; \; \text{ otherwise }. \\
\end{split}
\right.
\]
%

The most classical example of an isogeny is the  scalar multiplication by $n$ map  $[n] : \, \A \to \A$.  The kernel of $[n]$ is a group scheme of order $n^{2\dim \A}$.
Denote  by $\A [n]$  the group $\ker  [n] (\bar{k}) $.   The elements in $\A[n]$  are called $n$-\textit{torsion points} of $\A$.  Let $f : \A \to \B$ be a degree $n$ isogeny.  Then there exists an isogeny $\hat f : \B \to \A$ such that
\[ f \circ \hat f = \hat f \circ f = [n]. \]
Next we consider the case when $\chara k = p$.  Let $\A/k$ be an Abelian variety, $p = \chara k$, and $\dim \A= g$. 

\begin{itemize}
\item[i)] If $p \nmid \, n$, then $[n]$ is separable, $\# \A[n]= n^{2g}$ and $\A[n]\iso (\Z/n\Z)^{2g}$.

\item[ii)] If $p \mid n$, then $[n]$ is inseparable.  Moreover, there is an integer $0 \leq i \leq g$ such that 
\[ \A [p^m] \iso (\Z/p^m\Z)^i, \; \text{for all } \; m \geq 1. \]
\end{itemize}

If $i=g$ then $\A$ is called \textit{ordinary}.  If $\A[p^s](\bar k)= \Z/p^{ts}\Z$ then the abelian variety has \textit{$p$-rank} $t$. If $\dim \A=1$ (elliptic curve) then it is called \textit{supersingular} if it has $p$-rank 0.    An abelian variety $\A$ is called \textit{supersingular} if it is isogenous to a product of supersingular elliptic curves. 
If $\dim \A \leq 2$ and $\A$ has $p$-rank 0 then $\A$ is supersingular.   This is not true for $\dim \A\geq 3$.
\subsection{Jacobian varieties}
Let $\X$ be a curve of positive genus and assume that there exists a $k$-rational point  $P_0\in \X(k)$ with attached prime divisor $\p_0$.
There exists an abelian variety $\Jac_k (\X)$ defined over $k$ and a uniquely determined  embedding 
\[
\phi_{P_0}: \X \rightarrow \Jac_k (\X) \; \text{ with } \;  \phi_{P_0}(P_0) = 0_{\Jac_k (\X)} 
\]
such that
\begin{enumerate}
\item for  all extension fields  $L$ of $k$ we get $\Jac_L \X=\Pic^0_{\X_L}(L)$  where this equality is given in a functorial way and

\item if $\A$ is an Abelian variety and $\eta: \X \rightarrow \A$ is a morphism sending $P_0$ to $0_\A$ then there exists a uniquely determined  homomorphism 
$\psi:   \Jac (\X)     \rightarrow \A$ with $ \psi \circ \phi_{P_0}= \eta$.
\end{enumerate}
$\Jac (\X)$ is uniquely determined by these conditions and is called the \textit{Jacobian variety} of $\X$. The map $\phi_{P_0}$ is given by sending a prime divisor $\p$ of degree $1$ of $\X_L$ to the class of $\p-\p_0$ in $\Pic^0_{\X_L}(L)$.  

Let $L/k$ be a finite algebraic extension. Then the Jacobian variety $\Jac_L {\X}$ of $\X_L$ is the scalar extension of $\Jac \X$ with $L$, hence a fiber product with projection $p$ to $\Jac \X$. The norm map is $p_*$, and the conorm  map is $p^*$.   By universality we get that 
if  $f:  \X \rightarrow \mathcal \DD$ is a surjective morphism of curves sending $P_0$ to $Q_0$,  then there is a uniquely determined surjective homomorphism  $f_*:\Jac  \X     \rightarrow \Jac {\DD}$,   such that $f_*\circ \phi_{P_0}=\phi_{Q_0}$.
A useful observation is that 
if  $\Jac \X$ is a simple abelian variety, and $\eta: \X\rightarrow \DD$ is a   separable	cover of degree $>1$, then $\DD$ is the projective line.  For more details on the general setup see \cite{frey-shaska} among many other authors. 

\subsection{Jacobian surfaces}\label{sect-3}
Abelian varieties of dimension 2 are often called Abelian (algebraic) surfaces.  We focus on Abelian surfaces which are Jacobian varieties.  Let $\X$ be a genus 2 curve defined over a field $k$. Then its gonality is $\gamma_\X = 2$. Hence, genus 2 curves are hyperelliptic and we denote the hyperelliptic projection by  $\pi : \X \to \P^1$.   By the Hurwitz formula this covering has $r=6$ branch points which are images of the Weierstrass points of $\X$.  The moduli space has dimension $r-3=3$. 
  
The arithmetic of the moduli space of genus two curves was studied by Igusa in his seminal paper \cite{Ig} expanding on the work of Clebsch, Bolza, and others.  Arithmetic invariants $J_2, J_4, J_6, J_8, J_{10}$ determine uniquely the isomorphism class of a genus two curve. Two genus two curves $\X$ and $\X^\prime$ are isomorphic over $\bar k$  if and only if there exists $\l \in \bar k^\star$ such that $J_{2i} (\X) = \l^{2i} J_{2i} (\X^\prime)$, for $i=1,\dots, 5$.
If $\chara k \neq 2$ then the invariant $J_8$ is not needed. 

From now on we assume  $\chara k \neq 2$.   Then   $\X$ has an affine Weierstrass equation 
\begin{equation}\label{eq-g-2}
y^2=f(x)= a_6 x^6 + \dots + a_1 x + a_0, 
\end{equation}
over $\bar k$, with discriminant $\Delta_f = J_{10} \neq 0$.   The moduli space $\M_2$  of genus 2 curves, via the Torelli morphism, can be identified with the moduli space of the principally polarized abelian surfaces $\mathbb A_2$ which are not products of elliptic curves. Its compactification $\mathbb A_2^\star$ is the weighted projective space ${\mathbb{WP}}^3_{(2, 4, 6, 10)} (k)$  via the Igusa invariants $J_2, J_4, J_6,  J_{10}$.   Hence, 
\[ \mathbb A_2 \iso {\mathbb{WP}}^3_{(2, 4, 6, 10)} (k) \setminus \{ J_{10} = 0 \}.\]
Given a moduli point $\p \in \M_2$,   we can recover the equation of the corresponding curve over a minimal field of definition following   \cite{m-sh}. 


\def\AA{\mathcal A}

%
It is well known that a  map of algebraic curves $f : X \to Y$ induces maps between their Jacobians  $f^*: \Jac Y \to \Jac X$ and $f_*: \Jac  X \to \Jac Y$. When $f$ is maximal then $f^*$ is injective and $\ker (f_*)$ is connected, see \cite{jsc} for more details. 

Let $\X$ be a genus $2$ curve and  $\psi_1 : \X \longrightarrow  E_1$ be a   degree $n$ maximal   covering from   $\X$   to an elliptic curve $E_1$.   Then ${\psi^*}_1 : E_1 \to \Jac (\X)$ is injective and the kernel of $\psi_{1,*} : \Jac (\X) \to E_1$ is an elliptic curve which we denote by $E_2$. For a fixed Weierstrass point $P \in \X$, we can embed $\X$ to its Jacobian via
\begin{equation}
\begin{split}
i_P: \X & \longrightarrow  \Jac (\X) \\
 x & \to [(x)-(P)]
\end{split}
\end{equation}
Let $g : E_2 \to \Jac (\X)$ be the natural embedding of $E_2$ in $\Jac (\X)$, then there exists $g^*: \Jac (\X) \to E_2$. Define $\psi_2=g^*\circ i_P: \X \to E_2$. So we have the following exact sequence
\begin{equation}\label{red-Jacobians}
 0 \to E_2 \buildrel{g}\over\longrightarrow  \Jac (\X) \buildrel{\psi_{1,*}}\over\longrightarrow  E_1 \to 0. 
\end{equation} 
The dual sequence is also exact
\[ 0 \to E_1 \buildrel{\psi_1^*}\over\longrightarrow  \Jac (\X) \buildrel{g^*}\over\longrightarrow  E_2 \to 0. \]
If $\deg (\psi_1)=2$ or it is an odd number then the maximal covering $\psi_2: \X \to E_2$ is unique (up to isomorphism of elliptic curves).  
The Hurwitz space $\H_\sigma$ of such covers is embedded as a subvariety of the moduli space of genus two curves $\M_2$; see \cite{deg3} for details. It is a $2$-dimensional subvariety of $\M_2$ which we denote by $\L_n$. An explicit equation for $\L_n$, in terms of the arithmetic invariants of genus $2$ curves, can be found in \cite{deg2} or \cite{m-sh} for $n=2$, in \cite{deg3} for $n=3$, and in \cite{deg5} for $n=5$. 
From now on, we will say that a genus $2$ curve $\X$  has an $(n, n )$-decomposable Jacobian if $\X$ is as above and the elliptic curves $E_i$, $i=1, 2$ are called the components of $\Jac (\X)$.  

%
For every positive integer $D\equiv 0,1 \pmod 4$ there is a Humbert surface $H_D$ in $\M_2$ which parametrizes curves $\X$ whose Jacobians admit an optimal action on $\O_D$; see \cite{HM95}.   Points on $H_{n^2}$ parametrize curves whose Jacobian admits an $(n, n)$-isogeny to a product of two elliptic curves. Such curves are the main focus of our study. We have the following result; see \cite[Prop. 2.14]{lombardo}.

\begin{prop}  
$\Jac (\X)$ is a geometrically simple Abelian variety if and only if it is not $(n, n)$-decomposable for any $n>1$. 
\end{prop}


A point lying on the intersection of two Humbert surfaces $H_{m^2} \cap H_{n^2}$  with $n\neq m$ corresponds either to a simple abelian surface with quaternionic multiplication by an (automatically indefinite) quaternion algebra over $\Q$, or to the square of an elliptic curve. This is in particular true for points lying on Shimura curves.

%
We study pairs $(E_1, E_2)$ of elliptic components and try to determine their number (up to isomorphism over $\bar k$) when they are  isogenous of degree $N$, for an integer $N \geq 2$.  
We denote by  $\phi_N (x,y)$ the $N$-th modular polynomial. Two  elliptic curves with $j$-invariants $j_1$ and $j_2$ are $N$-isogenous if and only if $\phi_N (j_1,j_2)=0$.    The equation $\phi_N (x, y)=0$ is the canonical equation of the modular curve $X_0 (N)$.  The equations of $X_0 (N)$ are well known. 
%

\subsection{Kummer surface and Shioda-Inose surface}\label{kummer}


To the Jacobian variety $\Jac ( \X )$ one can naturally attach two $\KT$ surfaces, the Kummer surface and a double cover of it called the Shioda-Inose surface.   
Let   $\mathfrak{i}$ be the involution automorphism on  the Jacobian given by $\mathfrak i : \p  \to - \p$. The quotient    $\Jac (\X) / \{\mathbb{I}, \mathfrak{i}\}$, is a singular surface  with sixteen ordinary double points. Its minimal resolution is called the \textit{Kummer surface} and denoted by $\Kum (\Jac(\X))$. We refer to \cite{m-sh, m-sh-2} for further details.  

The Inose surface, denoted by  $\Y := \SI (\Jac (\X))$,  was originally constructed as a double cover of the Kummer surface. Shioda and Inose then showed that the following diagram of rational maps, called a Shioda-Inose structure, induces an isomorphism of integral Hodge structures on the transcendental lattices of $\Jac(\X)$ and $\Y$, see \cite{SI} for more details. 
\[
\begin{tikzcd}[column sep=scriptsize]
\Jac(\X) \arrow[dr, dashrightarrow, "\pi_0"]   
& & \Y \arrow[dl, dashrightarrow, "\pi_1"] \\
& \Kum(\Jac(\X))
\end{tikzcd}
\]
%

A $\KT$ surface $\Y$ has \textit{Shioda-Inose} structure if it admits an involution fixing the holomorphic two-form, such that the quotient is the Kummer surface $\Kum (\A)$  of a principally polarized abelian surface and the rational quotient map $p: \Y \to \Kum(\A)$ of degree two induces a Hodge isometry between the transcendental lattices $T(\Y) (2)^3$ and $T (\Kum(\A))$, see \cite{m-sh} for more details. 

 An elliptic surface   $\mathcal E (\overline k (t))$ fibered over $\mathbb P^1$ with section can be described by a Weierstrass equation of the form 
\[ y^2 + a_1(t) xy + a_3(t) y = x^3 +a_2(t)x^2 +a_4(t) x+a_6(t)\] 
and $a_i(t)$ rational functions. If we assume that the elliptic fibration has at least one singular fiber then the following question is fundamental in arithmetic geometry. Find generators for the Mordell-Weil group of this elliptic surface fibered over $\mathbb P^1$.   
 
A theorem of Shioda and Tate connects the Mordell-Weil group $\mathcal E (\overline k (t))$ with the Picard group of the N\'eron-Severi group of $\mathcal E$. Therefore, determining  the Mordell-Weil group is equivalent to finding the Picard group of the N\'eron-Severi lattice of $\KT$ surface.

 A surface is called an \textit{elliptic fibration} if it is a minimal elliptic surface over $\mathbb P^1$ with a distinguished section $S_0$. The complete list of possible singular fibers has been given by Kodaira \cite{kodaira}.  To each elliptic fibration $\pi: \X \to \mathbb P^1$ there is an associated Weierstrass model $\overline \pi: \overline \X \to \mathbb P^1$  with a corresponding distinguished section $\overline S_0$ obtained by contracting all fibers not meeting $S_0$.  The fibers of  $\overline \X$ are all irreducible whose singularities are all rational double points, and $\X$ is the minimal desingularization. If we choose some $t \in \X$ as a local affine coordinate on $\mathbb P^1$, we can present $\overline \X$ in the Weierstrass normal form 
 \[ Y^2 = 4X^3 -g_2(t) X - g_3(t),\]
 where $g_2(t)$ and $g_3(t)$ are polynomials of degree respectively $4$ and $6$ in $t$.  

\section{$(n, n)$ reducible Jacobian surfaces}\label{sect-4}

Genus 2 curves with $(n, n)$-decomposable Jacobians are the most studied type of genus $2$ curves due to work of  Jacobi,  Hermite, et al.  They provide examples of genus two curves with large Mordell-Weil rank of the Jacobian, many rational points, nice examples of descent \cite{ants}, etc. 
Such curves have received new attention lately due to   interest in their use in cryptographic applications and their suggested use on post-quantum crypto-systems and random self-reducibility of discrete logarithm problem; see \cite{costello}.  A detailed account of applications of such curves in cryptography is provided in \cite{frey-shaska}.

Let $\X$ be a genus $2$ curve defined over an algebraically closed field $k$, $\chara k =0$,  $K$ the function field of $\X$,  and $\psi_1 : \X \longrightarrow  E$  a  degree $n$ covering from $\X$  to an elliptic curve $E$; see \cite{jsc} for the basic definitions. The covering $\psi_1 : \X \longrightarrow  E$ is called a \textit{maximal covering} if it does not factor through a nontrivial isogeny. 
We call $E$  a \emph{degree $n$ elliptic subcover} of $\X$. Degree $n$ elliptic subcovers occur in pairs, say  
  $(E_1, E_2)$. It is well known that there is an isogeny of degree $n^2$ between the Jacobian $\Jac (\X)$   and the product $E_1 \times E_2$. Such a curve $\X$ is said to have  $(n, n)$-decomposable  (or $(n, n)$-split)  Jacobian. The focus of this paper is on isogenies among the elliptic curves $E_1$ and $E_2$. 

 The locus of genus $2$ curves  $\X$ with $(n, n)$-decomposable Jacobian is denoted by $\L_n$. When $n=2$ or $n$ an odd integer, $\L_n$  is a $2$-dimensional algebraic subvariety of the moduli space $\M_2$ of genus two curves; see \cite{jsc} for details.  Hence, we can get an explicit equation of $\L_n$ in terms of the Igusa invariants $J_2, J_4, J_6, J_{10}$; see \cite{deg2} for $\L_2$,    \cite{deg3} for $\L_3$, and  \cite{deg5} for $\L_5$.  There is a more recent paper on the subject \cite{kumar} where results of \cite{deg3, deg5} are confirmed and equations for   $n>5$ are studied.  
 

\subsection{$(2, 2)$ reducible Jacobian surfaces}  


Let $\X/k$ as above and $\F=k(\X)$ its function field. We assume that $k$ is algebraically closed and $\chara k \neq 2$.  Since degree $2$ coverings correspond to Galois extensions of function fields, the elliptic subcover is fixed by an involution in $\Aut (\F/k)$.  There is a group theoretic aspect of the $n=2$ case which was discussed in detail in \cite{deg2}. The number of elliptic subcovers in this case corresponds to the number of non-hyperelliptic involutions in $\Aut (\F/k)$, which are called  \textit{elliptic involutions}.  The equation of $\X$ is given by 
\[ 
Y^2=X^6-s_1X^4+s_2X^2-1
\]
see \cite{lubjana, beshaj-2} for minimal models of such equations.  Hence, for $(s_1,s_2)\in k^2$, such that the corresponding discriminant is nonzero,  we have a genus $2$ curve  $\X_{(s_1, s_2)}$ and two corresponding elliptic subcovers. Two such  curves  $(\X_{(s_1, s_2)},  \e_{s_1, s_2})$ and  $(\X_{(s_1', s_2')}, \e_{s_1', s_2'})$ are isomorphic if and only if their dihedral invariants $u$ and $v$ are the same; \cite{deg2}. Thus, the points $(s_1, s_2)\in k^2$ correspond to elliptic involutions of $\Aut \X$ while the points $(u, v)\in k^2$ correspond to elliptic involutions of  the reduced automorphism group $\bAut \X$. 


Let $\X$ be a genus $2$ curve, $\Aut (\X)$ its automorphism group, $\sigma_0$ the hyperelliptic involution, and    $\bAut (\X) := \Aut (\X)/ \< \sigma_0\>$ the reduced automorphism group. 
If $\Aut (\X)$ has another involution $\sigma_1$, then the quotient space $\X/\< \sigma_1\>$ has genus one. We call such an involution an \textit{elliptic involution}. There is another  elliptic involution   $\sigma_2:= \sigma_0\, \sigma_1$.    So the elliptic involutions come naturally in pairs. The corresponding coverings $\psi_i: \X \to \X/ \<\sigma_i\>$, $i=1, 2$,  are the maximal covers as above and  $E_i:=\X/ \<\sigma_i\>$ the elliptic  subcovers of $\X$ of degree $2$.  Also the corresponding Hurwitz space of such coverings is an irreducible algebraic variety which is embedded into $\M_2$.  We denote its image in $\M_2$ by $\L_2$.    The following was proved in \cite{deg2}. 


\begin{lem} \label{lem1} Let $\X$ be a genus $2$ curve and $\sigma_0$ its hyperelliptic involution. 
If $\sigma_1$ is an elliptic involution of $\X$, then so is  $\sigma_2=\sigma_1 \sigma_0$. Moreover,  $\X$ is isomorphic to a curve with affine equation
\begin{equation}\label{eq-1}
Y^2=X^6-s_1X^4+s_2X^2-1
\end{equation}
for some $s_1, s_2 \in k$ and  $\Delta :=27-18s_1s_2-s_1^2s_2^2+4s_1^3+4s_2^3\neq 0$. The equations for the elliptic subcovers  $E_i = \X/ \<\sigma_i\>$, for $i=1, 2$, are given by 
\[ E_1: \; y^2 = x^3-s_1x^2+s_2 x-1,  \;   \text{ and } \;     E_2:  \, y^2 = x \, (x^3-s_1x^2+s_2 x-1  )  \]
\end{lem}

In \cite{deg2} it was shown that   $\X$ is determined up to a  coordinate change by the subgroup $H\iso D_3$ of $SL_2 (k)$ generated by $\tau_1: X\to \e_6X$ and $\tau_2: X\to \frac 1 X$, where $\e_6$ is a primitive 6-th root of unity. Let $\e_3: = \e_6^2$. The coordinate change by $\tau_1$ replaces $s_1$ by $\e_3 s_1$ and $s_2$ by $\e_3^2 s_2$. The coordinate change by $\tau_2$ switches $s_1$ and $s_2$. Invariants of this $H$-action are:
\begin{equation}\label{eq2} 
u:=s_1 s_2, \quad v:=s_1^3+s_2^3 
\end{equation}
%
%
%
which are known in the literature as  \textit{dihedral invariants}.  The map 
\[
(s_1, s_2) \mapsto (u, v),
\] 
is a branched Galois covering with group $S_3$  of the set $\{ (u,v)\in k^2 : \Delta (u, v) \neq 0\}$ by the corresponding open subset of $(s_1, s_2)$-space  if $\chara (k) \ne 3$. In any case, it is true that if $s_1, s_2$ and $s_1', s_2'$ have the same $u,v$-invariants then they are conjugate under $\< \tau_1, \tau_2\>$.

If char$(k)=3$ then $u=u'$ and $v=v'$ implies $s_1^3 s_2^3 = s_1'^3 s_2'^3$ and $s_1^3+s_2^3= s_1'^3+s_2'^3$, hence 
$(s_1^3,s_2^3)= (s_1'^3,s_2'^3)$  or
$ (s_1^3,s_2^3)= (s_2'^3,s_1'^3)$. But this implies $(s_1,s_2)= (s_1',s_2')$
 or $(s_1,s_2)= (s_2',s_1')$.

For $(s_1,s_2)\in k^2$ with $\Delta\neq 0$, \cref{eq-1} defines a  genus 2 field  $\F_{s_1, s_2}= k(X,Y)$. Its reduced automorphism group contains the elliptic involution $\e_{s_1, s_2}: X \mapsto -X$. Two such  pairs $(\F_{s_1, s_2}, \e_{s_1, s_2})$ and  $(\F_{s_1', s_2'}, \e_{s_1', s_2'})$ are isomorphic if and only if $u=u'$ and $v=v'$ (where $u,v$ and $u',v'$ are associated with $s_1, s_2$ and $s_1', s_2'$, respectively, by \cref{eq2}).
However, the ordered pairs $(u, v)$ classify  the isomorphism classes of such elliptic subfields as it can be seen from the following theorem proved in \cite{deg2}.

\begin{prop}\label{sect1_thm}
i) The $(u,v)\in k^2$ with $\Delta\neq 0$ bijectively parameterize the isomorphism classes of pairs $(\F,\e)$ where $\F$ is a genus $2$ field and $\e$ an elliptic involution of $\bAut (\F)$. 

ii) The $(u,v)$ satisfying additionally
\begin{equation}\label{V_4}
(v^2-4u^3)(4v-u^2+110u-1125)\neq 0
\end{equation}
bijectively parameterize the isomorphism classes of genus $2$ fields with $\Aut (\F)\iso V_4$; equivalently, genus $2$ fields having exactly $2$ elliptic subfields of degree $2$. 
\end{prop}

Our goal  is to investigate when the pairs of elliptic subfields $\F_{s_1, s_2}$  (respectively isomorphism classes $(\F,\e)$) are isogenous. We want to find if that happens when $\X$ is defined over a number field $K$.  Hence, the following result is crucial. 

\begin{prop}\label{prop-deg2-invariants}

Let $K$ be a number field and $\X/K$ a genus $2$ curve with
\[
\operatorname{Aut}_{\overline K}(\X)\iso V_4.
\]
Then its dihedral invariants $u,v$ belong to $K$, and the $j$-invariants of its complementary elliptic components satisfy the quadratic in \cref{eq-quadratic}. Put
\begin{equation}\label{S-u-v}
S_2(u,v):=v^2-4u^3.
\end{equation}
The discriminant of this quadratic is
\[
\left(\frac{256(v-9u+27)}{\Delta}\right)^2 S_2(u,v).
\]
Consequently, both $j$-invariants belong to $K$ if and only if either $v-9u+27=0$ or $S_2(u,v)$ is a square in $K$.

\end{prop}

\proof   
The reduced elliptic involution is unique under the stated automorphism hypothesis, so the associated invariants $u,v$ are Galois invariant. Let $j_1,j_2$ denote the $j$-invariants of the elliptic components from \cref{lem1}. They are given in terms of $s_1,s_2$ by

\[
\begin{split}
j_1 & = -256\,{\frac { \left( {s_1}^{2}-3\,s_2 \right) ^{3}}{-{s_1}^{2}{s_2}^{2}+4\,{s_1}^{3}+4\,{s_2}^{3}-18\,s_1\,s_2+27}}\\
j_2 & = 256\,{\frac { \left( -{s_2}^{2}+3\,s_1 \right) ^{3}}{-{s_1}^{2}{s_2}^{2}+4\,{s_1}^{3}+4\,{s_2}^{3}-18\,s_1\,s_2+27}}\\
\end{split}
\]

Direct computation shows that they satisfy the quadratic

\begin{equation}\label{eq-quadratic}
j^2 - \left(  256    \frac{v^2 - 2u^3+54u^2-9uv-27v}{\Delta}  \right) j +  65536   \frac{(u^2+9u-3v)^3}{\Delta^2} =0
\end{equation}

where $\Delta=\Delta(u,v)=u^2-4v+18u-27$. Its discriminant has the asserted factorization. The final assertion follows from the quadratic discriminant criterion.

\qed

See \cite{main} for details, where an explicit equation of $\X$ is provided with coefficients as rational functions in  $u$ and  $v$,  or \cite{m-sh} for a more general setup. 
 
To determine whether the splitting itself descends, we use the two nonhyperelliptic involutions.

\begin{lem}\label{split-2}

Under the hypotheses of \cref{prop-deg2-invariants}, $\Jac(\X)$ admits a polarized $(2,2)$-splitting over $K$ if and only if $S_2(u,v)$ is a square in $K$.

\end{lem}

\proof

Over $\overline K$, choose the normal form
\[
y^2=x^6-s_1x^4+s_2x^2-1.
\]
Since the nontrivial reduced involution is unique, Galois transitions between these normal forms act on $x$ by a scaling or an inversion. A scaling preserves the two nonhyperelliptic involutions individually and preserves $s_1^3,s_2^3$ separately. An inversion interchanges both pairs.

Now $s_1^3,s_2^3$ are the distinct roots of
\[
T^2-vT+u^3=0.
\]
Thus the two nonhyperelliptic involutions are individually defined over $K$ precisely when $v^2-4u^3$ is a square in $K$. Their quotient curves and the induced splitting into the Jacobians of those quotients then descend to $K$. Conversely, a polarized $(2,2)$-splitting over $K$ gives a degree-two elliptic quotient structure defined over $K$, and hence one of these involutions.

\qed

We define the following surface 
\begin{equation} 
\S_2 : \; y^2 = S_2 (u, v), 
\end{equation}
where $S_2(u, v)$ is as in \cref{S-u-v}. 
%
%
Coefficients of \cref{eq-quadratic} can be expressed in terms of the Siegel modular forms or equivalently in terms of the Igusa arithmetic invariants; see \cite{thesis} or \cite{deg2}.  They were discovered independently in \cite{clingher-doran}, where  they are called \textit{modular invariants}. 

We consider the projection $\Phi:\S_2\to\L_2$ on the locus with geometric automorphism group $V_4$.

Then we have the following.

\begin{prop} 

On the locus with geometric automorphism group $V_4$, the projection
\[
\Phi:\S_2\longrightarrow\L_2,\qquad (u,v,y)\longmapsto(u,v),
\]
with $\S_2$ restricted to this locus, is a degree-two cover. Its two geometric points above $(u,v)$ correspond to the two choices of nonhyperelliptic involution. For $(u,v)\in K^2$ on this locus, the fiber consists of two $K$-rational points precisely when the splitting in \cref{split-2} is defined over $K$.

\end{prop}

\proof

The roots $s_1^3,s_2^3$ of $T^2-vT+u^3$ are distinct on this locus, and
\[
(s_1^3-s_2^3)^2=S_2(u,v).
\]
The scaling and inversion argument in the proof of \cref{split-2} identifies the two choices of sign with the two choices of nonhyperelliptic involution. The assertion about $K$-rational fibers follows from \cref{split-2}.

\qed

%
%

%
Next we turn our attention to isogenies between $E_1$ and $E_2$. 

\begin{prop}\label{prop-1}

Let $K$ be a number field and $N>1$. Let $G_N(u,v)$ be the numerator, in lowest terms and up to a nonzero scalar, of $\phi_N(j_1,j_2)$, expressed in $u,v$ using \cref{eq-quadratic}.

The $K$-rational points satisfying
\[
G_N(u,v)=0,\qquad
\Delta(u,v)(v^2-4u^3)(4v-u^2+110u-1125)\ne0
\]
correspond bijectively to geometric isomorphism classes of genus $2$ curves admitting a model over $K$, with geometric automorphism group $V_4$, whose complementary elliptic components are connected by a cyclic $N$-isogeny over $\overline K$.

\end{prop}

\proof

The symmetry of the modular polynomial expresses $\phi_N(j_1,j_2)$ in terms of $j_1+j_2$ and $j_1j_2$, which are given by \cref{eq-quadratic}. Its denominator is a power of $\Delta$, so on the stated locus its vanishing is equivalent to $G_N(u,v)=0$.

By the defining property of the modular polynomial, this vanishing is equivalent to the existence of a cyclic $N$-isogeny between the elliptic components over $\overline K$.

\Cref{sect1_thm} identifies the geometric isomorphism classes on the exact-$V_4$ locus with their dihedral invariants. \Cref{prop-deg2-invariants} shows that these invariants belong to $K$ for a curve defined over $K$. Conversely, a $K$-rational invariant pair determines a curve whose field of moduli is contained in $K$; the field-of-moduli result recalled above supplies a model over $K$. This proves the correspondence.

\qed

Let us now explicitly check whether elliptic components  of $\A$ are isogenous to each other. First we focus on the $d$-dimensional loci, for  $d\geq 1$. 

\begin{prop}\label{prop-2}

Let $K$ be a number field. For each $N=2,3,5,7$, there are infinitely many genus $2$ curves $\X$, up to geometric isomorphism, defined over $K$, with geometrically $(2,2)$-decomposable Jacobian and
\[
\operatorname{Aut}_{\overline K}(\X)\iso V_4,
\]
whose complementary elliptic components $E_1,E_2$ are $N$-isogenous over $\overline K$.

\end{prop}

\proof
Let us now check if elliptic components are isogenous   for $N=2, 3, 5, 7$. By replacing $j_1, j_2$ in the modular curve we get a curve   $F(s_1, s_2) =0 $.
This curve is symmetric in $s_1$ and $s_2$ and fixed by the $H$-action described in the preliminaries. Therefore, such a curve can be written in terms of $u$ and $v$, 
$ G_N (u, v) =0$.
 We display all the computations below. 

Let $N=2$.  $G_2 (u, v) $ is 
\[
G_2(u, v)= f_1 (u, v) \cdot f_2 (u, v)
\]
where $f_1$ and $f_2$ are
\begin{small}
\begin{equation}
\begin{split}
f_1 &
=-16v^3-81216v^2-892296v-2460375+3312uv^2+707616vu+3805380u+\\
&  18360vu^2 -1296162u^2 -1744u^3v-140076u^3+801u^4+256u^5 \\
\end{split}
\end{equation}
\end{small}
\begin{small}
\begin{equation}
\begin{split}
f_2 & =4096u^7+256016u^6-45824u^5v+4736016u^5-2126736vu^4+23158143u^4\\
&  -25451712u^3v-119745540u^3+5291136v^2u^2-48166488vu^2-2390500350u^2\\
& -179712uv^3+35831808uv^2+1113270480vu+9300217500u-4036608v^3\\
& -1791153000v-8303765625-1024v^4+163840u^3v^2-122250384v^2+256u^2v^3 \\
\end{split}
\end{equation}
\end{small}

For $N=2$, both components $f_1(u,v)=0$ and $f_2(u,v)=0$ are absolutely irreducible and have geometric genus zero.

Let $N=3$.   Then, from \cref{V_4}   and $\phi_3(j_1,j_2)=0$ we have:
\begin{equation}\label{iso_3}
(4v-u^2+110u-1125)\cdot g_1(u,v)\cdot g_2(u,v)=0  
\end{equation}
where $g_1$ and $g_2$ are
\begin{small}
\begin{equation}
\begin{split}
g_1 & =-27008u^6+256u^7-2432u^5v+v^4+7296u^3v^2-6692v^3u-1755067500u\\
 & +2419308v^3-34553439u^4+127753092vu^2+16274844vu^3-1720730u^2v^2\\
 & -1941120u^5+381631500v+1018668150u^2-116158860u^3+52621974v^2\\
 & +387712u^4v -483963660vu-33416676v^2u+922640625 \\
\end{split}
\end{equation}
\end{small}
\begin{small}
\begin{equation}
\begin{split}
g_2 & =291350448u^6-v^4u^2-998848u^6v-3456u^7v+4749840u^4v^2+17032u^5v^2\\
 &  +4v^5+80368u^8+256u^9+6848224u^7-10535040v^3u^2-35872v^3u^3+26478v^4u\\
 &  -77908736u^5v+9516699v^4+307234984u^3v^2-419583744v^3u-826436736v^3\\
 &  +27502903296u^4+28808773632vu^2-23429955456vu^3+5455334016u^2v^2\\
 &  -41278242816v+82556485632u^2-108737593344u^3-12123095040v^2\\
 & +41278242816vu+3503554560v^2u+5341019904u^5-2454612480u^4v\\
\end{split}
\end{equation}
\end{small}

For $N=3$, the factor $4v-u^2+110u-1125$ is excluded by the condition $\Aut(\X)\iso V_4$. The remaining components $g_1(u,v)=0$ and $g_2(u,v)=0$ are absolutely irreducible and have geometric genus zero.

For $N=5,7$, the modular pullback has three absolutely irreducible components in each case. Their degrees and geometric genera, together with those of the preceding cases after imposing the $V_4$ condition, are:
\[
\begin{array}{c|c|c}
N&\text{Degrees of the plane closures}&\text{Geometric genera}\\ \hline
2&5,7&0,0\\
3&7,9&0,0\\
5&4,11,15&0,0,0\\
7&6,15,21&0,0,0
\end{array}
\]
These results follow from exact factorization and geometric-genus computations. Each component also has a smooth rational point, as can be checked by restricting its equation to $v=9u-27$. Consequently, its smooth projective normalization is isomorphic over $\mathbb Q$ to $\mathbb P^1$.

The nonsingularity and automorphism conditions exclude only finitely many points from each retained component: none is contained in
\[
\Delta(u,v)(v^2-4u^3)(4v-u^2+110u-1125)=0.
\]
Thus each retained component contains infinitely many rational moduli points satisfying these conditions. By the field-of-moduli result recalled above for genus $2$ curves with extra involutions, the corresponding curves admit models over $\mathbb Q$, and hence over $K$. Since the dihedral invariants classify geometric isomorphism classes on the exact-$V_4$ locus, these give infinitely many distinct classes.

\qed

Next we consider the case when $|\Aut (\X) | > 4$. First notice that the invariants $j_1$ and $j_2$ are roots of the quadratic \cref{eq-quadratic}. 

Assume that $G=\operatorname{Aut}_{\overline K}(\X)\iso D_4$. The nonhyperelliptic involutions fall into two conjugacy classes. Write $E_1,E_2$ for the elliptic quotients arising from one class and $E_1^\prime,E_2^\prime$ for those arising from the other. Within each pair the quotients are isomorphic; the isomorphism classes represented by the two pairs may also coincide. Over $\overline K$, we may choose the even normal form with $s_1=s_2=t$, so that
\[
u=t^2,\qquad v=2t^3.
\]

%

%
 
In this case $\bG\iso V_4$, and the set of Weierstrass points is

\[ \mathcal W=\{\pm 1, \pm \sqrt{a}, \pm \sqrt{b}\}.\]

Here $b=a^{-1}$ and $t=a+a^{-1}+1$. The reduced involutions are $\tau_1:X\mapsto-X$, $\tau_2:X\mapsto1/X$, and $\tau_3:X\mapsto-1/X$. Since $\tau_1$ and $\tau_3$ fix no points of $\mathcal W$, they lift to elliptic involutions. Each determines a pair of isomorphic elliptic quotients. The $j$-invariant of the quotients associated with $\tau_1$ is the double root of \cref{eq-quadratic}, namely

\begin{equation}\label{D4-j1}
 j=256\frac{t^3}{t+1}. 
\end{equation}

To find the $j$-invariant of the quotients associated with $\tau_3$, consider the following degree $2$ covering.
The covering $\phi:\bP^1\to\bP^1$ satisfies $\phi(\pm1)=0$,
\[
\phi(\sqrt a)=\phi(-1/\sqrt a)=1,\qquad
\phi(-\sqrt a)=\phi(1/\sqrt a)=-1,
\]
and $\phi(0)=\phi(\infty)=\infty$. It is given by
\[
\phi(X)=\frac{\sqrt a}{a-1}\frac{X^2-1}{X}.
\]
Its critical values are
\[
q_1=-q_2=\frac{2i\sqrt a}{a-1}.
\]
For each $i=1,2$, the elliptic double cover corresponding to $E_i^\prime$ has branch values $\{0,1,-1,q_i\}$. Computing its $j$-invariant gives

\begin{equation}\label{D4-j2} 
j^\prime=-16\frac{(t-15)^3}{(t+1)^2}. 
\end{equation}

Then, we have the following result.

\begin{prop}\label{prop-3}

Let $K$ be a number field and $\X/K$ a genus $2$ curve with
\[
\operatorname{Aut}_{\overline K}(\X)\iso D_4.
\]
Let $E_i,E_i^\prime$, $i=1,2$, be the elliptic quotients arising from the two conjugacy classes of nonhyperelliptic involutions, as above. Then $E_i$ and $E_i^\prime$ are $2$-isogenous over $\overline K$. For each $N=3,5,7$, there are only finitely many such curves $\X$, up to geometric isomorphism, admitting a model over $K$, for which $E_i$ and $E_i^\prime$ are also $N$-isogenous over $\overline K$.

\end{prop}

\proof

Over $\overline K$, choose the normal form with $s_1=s_2=t$. Its discriminant satisfies
\[
\Delta=(t-3)^3(t+1),
\]
so $t\ne-1,3$. Using \cref{D4-j1,D4-j2}, direct substitution gives
\[
\phi_2(j,j^\prime)=0
\]
and

\begin{small}
\[
\begin{aligned}
\phi_3(j,j^\prime)={}&\frac{2^{16}}{(t+1)^9}
(t^2+138t+153)(t+5)^2(t^2-70t-55)^2\\
&\cdot(256t^4+240t^3+191745t^2+371250t+245025)\\
&\cdot(4096t^6-17920t^5+55909200t^4-188595375t^3\\
&\hspace{26mm}-4518125t^2+769621875t+546390625).
\end{aligned}
\]
\end{small}

The first identity proves the $2$-isogeny assertion.

For $N=3,5,7$, the reduced numerators of $\phi_N(j,j^\prime)$ are nonzero polynomials in $t$, of degrees $18,30,42$, respectively. Each therefore has only finitely many roots. Restricting to parameters defining smooth curves with geometric automorphism group exactly $D_4$ leaves finitely many admissible values. Each such value determines a geometric isomorphism class of genus $2$ curves, proving the assertion.

 This completes the proof.
\qed

\subsection{$(3, 3)$ reducible Jacobian surfaces}\label{sect-5}

In this section we focus on genus $2$ curves with $(3, 3)$-split Jacobians.   This case was studied in detail in \cite{deg3}, where it was  proved that if $\F$ is a genus $2$ field over $k$  and $e_3(\F)$ the number of $\Aut (\F/k)$-classes   of elliptic subfields of $\F$ of degree $3$,   then

i)  $e_3 (\F) =0, 1, 2$, or  $4$

ii)    $e_3(\F) \geq 1$ if and only if  the classical invariants of $\F$ satisfy  the irreducible equation   $f(J_2, J_4, J_6, J_{10})=0$ displayed  in  \cite[Appendix A]{deg3}.

There are exactly two genus $2$ curves (up to isomorphism) with $e_3(\F)=4$. The case $e_3(\F)=1$ (resp., $2$) occurs for a $1$-dimensional (resp., $2$-dimensional)  family of genus $2$ curves.  We are interested in the $2$-dimensional family, since the case $e_3 (\F)= 1$ is the singular locus of the case $e_3 (\F)=2$.

We let $\X$ be a genus 2 curve defined over $k=\bar k$, $\chara k \neq 2, 3$, and $\F:=k(\X)$ its function field. 
\begin{defn}
A {\bf non-degenerate pair} (resp., {\bf degenerate pair}) is a
pair $(\X,\E)$ such that $\X$ is a genus $2$ curve with a degree $3$
elliptic subcover $\E$ where $\psi: \X \to \E$ is ramified in two
(resp., one) places.  Two such pairs $(\X, \E)$ and $(\X' , \E')$
are called isomorphic if there is a $k$-isomorphism $\X \to \X'$
mapping $\E\to \E'$.
\end{defn}

On the following open chart, a non-degenerate pair $(\X,\E)$ is represented by

\begin{equation}\label{eq_F1_F2}
Y^2=(\v^2 X^3+\u \v X^2 +\v X+1)\, (4\v^2 X^3 +\v^2 X^2+2\v X+1),
\end{equation}
where $\u, \v  \in k$ and the discriminant
\[ \Delta = -16\, \v^{17} \, (\v -27) \, (27 \v + 4 \v^2 - \u^2 \v + 4\u^3 -18 \u \v )^3 \]
of the sextic is nonzero. We let
$R:=(27 \v + 4 \v^2 - \u^2 \v + 4\u^3 -18 \u \v )\neq 0$.
For $4\u - \v -9\neq 0$ the degree $3$ coverings are given by
$\phi_1(X,Y)\to (U_1, V_1)$ and $\phi_2(X,Y)\to (U_2, V_2)$ where
\begin{small}
\begin{equation}\label{covers}
\begin{split}
 U_1 & = \frac {\v X^2} {\v^2 X^3 + \u \v X^2 + \v X +1}, \\
 U_2  & = \frac {(\v X +3)^2\, \,  (\v  (4\u-\v-9) X +3\u -\v)} {\v \,(4\u -\v -9) (4 \v^2 X^3 + \v^2 X^2 + 2 \v X + 1) }, \\
   V_1 & =  Y \, \frac {\v^2 X^3 - \v X -2} {(\v^2 X^3 + \u \v X^2 + \v X +1)^2}, \\
  V_2  & =   (27- \v)^{ \frac 3 2} \, Y
\frac {\v^2 (\v-4\u+8)X^3 +\v (\v-4\u)X^2 -\v X-1  }
{(4 \v^2 X^3 + \v^2 X^2 + 2 \v X + 1)^2 }\\
\end{split}
\end{equation}
\end{small}
and the elliptic curves have equations:
\begin{small}
\begin{equation}\label{ell}
\begin{split}
\E_1: & \quad V_1^2= -\big(R \, U_1^3 - (12 \u^2 - 2 \u \v - 18 \v )U_1^2 +
(12\u -\v) U_1 -4\big) \\
\E_2: & \quad V_2^2=c_3 U_2^3 +c_2 U_2^2 + c_1 U_2 +c_0
\end{split}
\end{equation}
\end{small}
where
\begin{small}
\begin{equation}
\begin{split}
c_0 & = - (9\u -2 \v -27 )^3\\
c_1 & = (4\u -\v -9) \, (729 \u^2 + 54 \u^2 \v -972\u \v - 18\u \v^2 +189\v^2 +
729 \v +\v^3) \\
c_2 & =- \v \,  (4\u -\v -9)^2 \,  (54\u +\u \v -27\v)     \\
c_3 & =\v^2 \, (4\u -\v -9)^3 \\
\end{split}
\end{equation}
\end{small}

For a subfield $K$ containing $\u,\v$, the curve $\E_1$ and the covering $\phi_1$ are defined over $K$. Omitting the factor $(27-\v)^{3/2}$ in \cref{covers} gives a $K$-defined map to the quadratic twist
\[
W_2^2=\frac{c_3U_2^3+c_2U_2^2+c_1U_2+c_0}{(27-\v)^3}
\]
of $\E_2$, which has the same $j$-invariant.

The parameter map from \cref{eq_F1_F2} to $\L_3$ has generic degree $2$. On the chart $J_2\ne0$, it is written in absolute invariants as
\[
(\u,\v)\longmapsto(i_1,i_2,i_3).
\]
%



%

We define the following invariants of two cubic polynomials. For $F(X)=a_3X^3+a_2 X^2 + a_1 X+a_0$ and $G(X)=b_3 X^3+b_2X^2+b_1X+b_0$ define
\[ H(F,G) :=a_3 b_0 - \frac  1 3 a_2 b_1 + \frac 1 3 a_1 b_2 -a_0 b_3 \]
We denote by $R(F,G)$ the resultant of $F$ and $G$ and by $D(F)$ the discriminant of $F$. Also,
$r_1(F,G) =\frac  {H(F,G)^3} {R(F,G)}$, $r_2(F,G)=\frac {H(F,G)^4} {D(F)\, D(G)}$, and $ r_3=\frac {H(F,G)^2} {J_2(F\,G)}$. 

Invariants $r_1,r_2,r_3$ form a complete system of invariants for unordered pairs of cubics on the generic locus. For the two cubic factors $F,G$ of \cref{eq_F1_F2}, on the chart $\v(\v-27)R(\v-2\u-9)\ne0$, we use the normalization

\[
\begin{split}
\chi:= -r_1 &= 3^{-3} \cdot \frac {{\v}({\v}-9-2{\u})^3} {4{\v}^2-18{\u}{\v}+27{\v}-{\u}^2{\v}+4{\u}^3}\\
\psi := r_2 & = -2^{-4} \cdot 3^{-4} \frac  {{\v}({\v}-9-2{\u})^4}   {({\v}-27)(4{\v}^2-18{\u}{\v}+27{\v}-{\u}^2{\v}+4{\u}^3)},
\end{split}
\]

It was shown in \cite{deg3} that the function field of the locus $\L_3$, genus 2 curves with $(3,3)$ reducible Jacobians, is exactly $k(\chi, \psi)$. 

\begin{lem}  $k(\L_3) = k (\chi, \psi)$. 
These coordinates are used on the preceding generic chart.

\end{lem}

By eliminating $\u$ and $\v$ we have rational expressions of absolute invariants $i_1, i_2, i_3$ in terms of $\chi$ and $\psi$ as in  \cite[Eq.~(19)]{deg3}.
 
On the chart $J_2\ne0$, we can take the following weighted representative:

\[ \left[ J_2 : J_4 : J_6 :J_{10}   \right] = \left[ 1 : \frac 1 {144} i_1   :  \left( \frac 1 {5184} i_2   + \frac 1 {432} i_1 \right)   : \frac 1 {486} i_3   \right]
\]

With the above normalization of $\chi,\psi$, a weighted Igusa representative is given by

\begin{small}
\[
\begin{split}
J_2 = & \chi\, \left( {\chi}^{2}+96\,\chi\,\psi-1152\,{\psi}^{2} \right) \\
J_4 = &  \frac {\chi}  {2^6} \, \left( {\chi}^{5}+192\,{\chi}^{4}\psi+13824\,{\chi}^{3}{\psi}^{2}+ 442368\,{\chi}^{2}{\psi}^{3}+5308416\,\chi\,{\psi}^{4} \right.\\
& \left. +786432\,\chi\,{ \psi}^{3} +9437184\,{\psi}^{4} \right)    \\
J_6 = & \frac {\chi} {2^9 }  \left( 3\,{\chi}^{8}+864\,{\chi}^{7}\psi+94464\,{\chi}^{6}{\psi}^{2}+4866048\,{\chi}^{5}{\psi}^{3}+111476736\,{\chi}^{4}{\psi}^{4} \right. \\
& +509607936\,{\chi}^{3}{\psi}^{5}  -12230590464\,{\chi}^{2}{\psi}^{6}+ 1310720\,{\chi}^{4}{\psi}^{3}+155713536\,{\chi}^{3}{\psi}^{4} \\
& \left. -1358954496\,{\chi}^{2}{\psi}^{5}  -18119393280\,\chi\,{\psi}^{6}+4831838208\,{\psi}^{6} \right)  \\
J_{10} = & - 2^{30} \chi^3 \psi^9  \\
\end{split}
\]
\end{small}
It would be an interesting problem to determine for what values of $\chi$ and $\psi$ the curve $\X$ is defined over the field of moduli. 
\subsubsection{Elliptic components}

We express the j-invariants $j_i$ of the elliptic components $\E_i$ of $\A$, from \cref{ell}, in terms of $\u$ and $\v$ as follows:

\begin{small}
\begin{equation}\label{j_1}
\begin{split}
j_1 &= 16\v \frac {(\v\u^2+216\u^2-126\v\u-972\u+12\v^2+405\v)^3}   { (\v-27)^3(4\v^2+27\v+4\u^3-18\v\u-\v\u^2)^2}\\
j_2 &= -256 \frac {(\u^2-3\v)^3}{\v  (4\v^2+27\v+4\u^3-18\v\u-\v\u^2)} \\
\end{split}
\end{equation}
\end{small}
where $\v\neq 0, 27$
and $R\ne0$. Moreover, we can express $s=j_1+j_2$ and $t=j_1j_2$ in terms of the $\chi$ and $\psi$ invariants as follows:

\begin{lem}
 
On this chart, the $j$-invariants of the complementary elliptic components satisfy the following quadratic over $k(\chi,\psi)$:

\begin{equation}\label{eq_j_new}
j^2- s \, j+ t =0
\end{equation}
where 
\begin{small}

\begin{equation}\label{eq-j-deg3}
\begin{split}
s = & - \frac 1 {2^{24}\, \chi^4 \psi^6} \Big( \chi^{12} +480\,\chi^{11}\,\psi +101376\,\chi^{10}\,\psi^{2} \\
& +12386304\,\chi^{9}\,\psi^{3} +966131712\,\chi^{8}\,\psi^{4} +49941577728\,\chi^{7}\,\psi^{5} \\
& +1712282664960\,\chi^{6}\,\psi^{6} +37572373905408\,\chi^{5}\,\psi^{7} +479047767293952\,\chi^{4}\,\psi^{8} \\
& +2705210921189376\,\chi^{3}\,\psi^{9} -258048\,\chi^{8}\,\psi^{3} -38928384\,\chi^{7}\,\psi^{4} \\
& +1528823808\,\chi^{6}\,\psi^{5} +619683250176\,\chi^{5}\,\psi^{6} +45595641249792\,\chi^{4}\,\psi^{7} \\
& +1408964021452800\,\chi^{3}\,\psi^{8} +16231265527136256\,\chi^{2}\,\psi^{9} +7247757312\,\chi^{4}\,\psi^{6} \\
& +7827577896960\,\chi^{3}\,\psi^{7} +901736973729792\,\chi^{2}\,\psi^{8} +32462531054272512\,\chi\,\psi^{9} \\
& +21641687369515008\,\psi^{9} \Big) \\
t = & \frac 1 {2^{36}\, \chi^3 \psi^9} \Big( \chi^{5} +192\,\chi^{4}\,\psi +13824\,\chi^{3}\,\psi^{2}  +442368\,\chi^{2}\,\psi^{3} +5308416\,\chi\,\psi^{4} \\
	& -1179648\,\chi\,\psi^{3}  -84934656\,\psi^{4} \Big)^3
\end{split}
\end{equation}

\end{small}
\end{lem}

\proof Substitute $j_1$ and $j_2$ as in \cref{j_1} in \cref{eq_j_new}. \endproof

\begin{rem}
On this birational chart, a model of $\X$ over $K$ gives $\chi,\psi\in K$. A rational moduli point does not by itself supply a model of $\X$ over $K$.
\end{rem}

Invariants $s$ and $t$ are \textit{modular invariants} similar to the $n=2$ case and can be expressed in terms of the Siegel modular forms or equivalently in terms of the Igusa arithmetic invariants.

Let $K$ be a perfect field of characteristic different from $2,3$, and let $(\chi,\psi)\in K^2$ lie on the preceding smooth moduli chart, with $\chi\psi\ne0$.

 The discriminant of the quadratic in \cref{eq_j_new} is given by
\begin{small}

\begin{equation} 
\begin{split}
\Delta (\chi, \psi) = & \frac 1 {2^{48} \, \chi^{8} \psi^{12}} \Big( \chi^{8} +288\,\chi^{7}\,\psi +34560\,\chi^{6}\,\psi^{2} \\
& +2211840\,\chi^{5}\,\psi^{3} +79626240\,\chi^{4}\,\psi^{4} +1528823808\,\chi^{3}\,\psi^{5} \\
& +12230590464\,\chi^{2}\,\psi^{6} +110592\,\chi^{4}\,\psi^{3} +31850496\,\chi^{3}\,\psi^{4} \\
& +2293235712\,\chi^{2}\,\psi^{5} +48922361856\,\chi\,\psi^{6} +48922361856\,\psi^{6} \Big)^2 \\
& \cdot \Big( \chi^{8} +384\,\chi^{7}\,\psi +59904\,\chi^{6}\,\psi^{2} \\
& +4866048\,\chi^{5}\,\psi^{3} +217645056\,\chi^{4}\,\psi^{4} +5096079360\,\chi^{3}\,\psi^{5} \\
& +48922361856\,\chi^{2}\,\psi^{6} -753664\,\chi^{4}\,\psi^{3} -113246208\,\chi^{3}\,\psi^{4} \\
& -2038431744\,\chi^{2}\,\psi^{5} +195689447424\,\chi\,\psi^{6} +195689447424\,\psi^{6} \Big)
\end{split}
\end{equation}

\end{small}

Let $A_3(\chi,\psi)$ denote the polynomial that is squared in this display. The second factor will be denoted by $S_3(\chi,\psi)$, so
\[
\Delta(\chi,\psi)=\left(\frac{A_3(\chi,\psi)}{2^{24}\chi^4\psi^6}\right)^2S_3(\chi,\psi).
\]
If $A_3(\chi,\psi)\ne0$, the discriminant is a square precisely when $S_3(\chi,\psi)$ is a square in $K$. At $A_3(\chi,\psi)=0$ the discriminant is zero regardless of $S_3(\chi,\psi)$. As in the case $n=2$, define the surface

\begin{equation} 
\S_3 : \qquad y^2=S_3 (\chi, \psi), 
\end{equation}
where
\begin{equation}
\begin{split}
S_3 (\chi, \psi) :=  &\chi^{8} +384\,\chi^{7}\,\psi +59904\,\chi^{6}\,\psi^{2}  +4866048\,\chi^{5}\,\psi^{3} +217645056\,\chi^{4}\,\psi^{4} \\			
			& +5096079360\,\chi^{3}\,\psi^{5}  +48922361856\,\chi^{2}\,\psi^{6} -753664\,\chi^{4}\,\psi^{3} \\
			&	-113246208\,\chi^{3}\,\psi^{4}  -2038431744\,\chi^{2}\,\psi^{5} +195689447424\,\chi\,\psi^{6} \\
			&	+195689447424\,\psi^{6}
\end{split}
\end{equation}

is the second factor in the discriminant. Projection gives a generically degree-two cover of the chosen birational chart of $\L_3$:
\[
\Phi:\S_3\longrightarrow\L_3,\qquad(\chi,\psi,y)\longmapsto(\chi,\psi).
\]
Over $A_3(\chi,\psi)S_3(\chi,\psi)\ne0$, its two geometric points label the two orderings of the complementary $j$-pair, since
\[
j_\pm=\frac{s}{2}\pm\frac{A_3(\chi,\psi)\,y}{2^{25}\chi^4\psi^6}.
\]
The fiber above the specified $(\chi,\psi)$ has a $K$-rational point precisely when $S_3(\chi,\psi)$ is a square in $K$. At $S_3(\chi,\psi)=0$ the two points coalesce. At $A_3(\chi,\psi)=0$ the $j$-values already belong to $K$, even if that fiber has no $K$-rational point.

\begin{lem}\label{split-3}

Let $K$ be a perfect field of characteristic different from $2,3$, and let $(\chi,\psi)\in K^2$ be a point of the preceding smooth moduli chart, with
$\chi\psi\ne0$. The two complementary elliptic $j$-invariants belong to $K$ if
and only if either $A_3(\chi,\psi)=0$ or $S_3(\chi,\psi)$ is a square in $K$.
Equivalently, this is the condition that their geometric elliptic
isomorphism classes admit models over $K$.

Suppose, in addition, that the data are represented by a curve $\X/K$
and that its chosen unordered pair of complementary degree-three elliptic
subvarieties is Galois invariant, with distinct $j$-invariants. Then the
associated polarized $(3,3)$-splitting of $\Jac(\X)$ is defined over $K$
if and only if $S_3(\chi,\psi)$ is a square in $K$.

\end{lem}

\proof  

The first assertion follows from the discriminant criterion for the
quadratic in \cref{eq_j_new}. If $A_3(\chi,\psi)=0$, its double root is $s/2\in K$.
Otherwise its discriminant is a square precisely when $S_3(\chi,\psi)$ is a
square. An elliptic geometric isomorphism class admits a model over $K$
precisely when its $j$-invariant belongs to $K$.

For the second assertion, Galois can only permute the two chosen elliptic
subvarieties. Since their $j$-invariants are distinct, both are individually
Galois invariant exactly when both $j$-invariants belong to $K$. The
subvarieties and their induced polarized isogeny then descend to $K$.
Conversely, a splitting over $K$ gives both elliptic $j$-invariants in $K$.
Distinctness excludes the zero discriminant, so $A_3(\chi,\psi)\ne0$ and the
first assertion gives the stated square criterion.

\qed

The unordered-pair hypothesis is automatic when the complementary pair is unique. The equal-$j$ case requires separate analysis. The second assertion concerns the polarized splitting of the Jacobian; descent of particular maps to pointed elliptic curves can additionally involve genus-one target torsors.

A rational moduli point need not give a genus $2$ model over its field of moduli. The reconstruction question for $\X$ is separate from the splitting assertion above; see \cite{m-sh}.

\subsubsection{Isogenies between the elliptic components}
Now let us consider the case when $n=3$. In an analogous way with the case $n =2$ we will study the locus $\phi_N (x, y)=0$ which represents the modular curve $X_0 (N)$. For $N$ prime, two  elliptic curves $E_1,~E_2$  are $N$-isogenous if and only if  $\phi_N (j(E_1),j(E_2))=0$.  We will consider the case when $N= 2, 3, 5$, and $7$.  We will omit part of the formulas since they are too big to display.  

\begin{prop}\label{prop-4}

Let $K$ be a number field. There are infinitely many genus $2$ curves, up to geometric isomorphism, admitting models over $K$, with geometrically $(3,3)$-split Jacobians whose complementary elliptic components are $5$-isogenous over $\overline K$.

\end{prop}

\proof  Let $\phi_5 (x, y)$ be the modular polynomial of level $5$.  As in the previous section, we let $s=x+y$ and $t=xy$.  Then, $\phi_5 (x, y)$ can be written in terms of $s, t$. We replace $s$ and $t$ by expressions in \cref{eq-j-deg3}.  

We get a curve in $\chi$, $\psi$ with three irreducible components over $\Q$, of degrees $4, 24, 56$. The component of degree $4$ is
\[ \chi^2 (\chi + 48 \psi)^2 = 2^{14} \psi^3. \]
Its preimage in the $(\u, \v)$-plane is the rational curve
\[ \u = \frac {4 (2m-1)(m+1)} {m^2}, \qquad \v = \frac {16 (m+1)^3} {m (m^2+4)}, \]
for which
\[ j_1 = - \frac {4096\, (31 m^2 + 57 m + 1)^3} {m\, (11 m - 4)^5}, \qquad j_2 = - \frac {4096\, (m^2 - 3 m + 1)^3} {m^5 (11 m - 4)} \]
and $\phi_5 (j_1, j_2)=0$. The discriminant of the sextic in \cref{eq_F1_F2} is a nonzero constant multiple of $m^{-36} (11m-4)^4 (m-2)^{20} (m+1)^{60} (m^2+4)^{-24}$. Hence, for every $m \in K\setminus \{0, -1, 2, 4/11\}$ with $m^2+4\neq 0$, \cref{eq_F1_F2} is a genus $2$ curve defined over $K$ whose elliptic components are $5$-isogenous. The rational function $j_2(m)$ is nonconstant, and a fixed genus $2$ curve has only finitely many degree-three elliptic subcovers; hence the family contains infinitely many geometric isomorphism classes. This completes the proof.  
\qed

Let us now consider the other cases. 
 
Let $G_N^{(3)}(\chi,\psi)$ be the reduced numerator of $\phi_N(j_1,j_2)$ after using \cref{eq-j-deg3}. For $N=2$, factorization over $\mathbb Q$ gives

\[ g_1 (\chi, \psi) \cdot g_2(\chi, \psi) =0, \]
where $g_1 (\chi, \psi)=0$ is a genus zero component given by 
\begin{small}

\begin{equation}
\begin{split}
& \chi^{12} +432\,\chi^{11}\,\psi +82944\,\chi^{10}\,\psi^{2} \\
& +9289728\,\chi^{9}\,\psi^{3} +668860416\,\chi^{8}\,\psi^{4} +32105299968\,\chi^{7}\,\psi^{5} \\
& +1027369598976\,\chi^{6}\,\psi^{6} +21134460321792\,\chi^{5}\,\psi^{7} +253613523861504\,\chi^{4}\,\psi^{8} \\
& +1352605460594688\,\chi^{3}\,\psi^{9} -2322432\,\chi^{8}\,\psi^{3} -525533184\,\chi^{7}\,\psi^{4} \\
& -45864714240\,\chi^{6}\,\psi^{5} -1834588569600\,\chi^{5}\,\psi^{6} -26418075402240\,\chi^{4}\,\psi^{7} \\
& +253613523861504\,\chi^{3}\,\psi^{8} +8115632763568128\,\chi^{2}\,\psi^{9} +2348273369088\,\chi^{4}\,\psi^{6} \\
& +190210142896128\,\chi^{3}\,\psi^{7} +4057816381784064\,\chi^{2}\,\psi^{8} +16231265527136256\,\chi\,\psi^{9} \\
& +10820843684757504\,\psi^{9} =0,
\end{split}
\end{equation}

\end{small}

The cofactor $g_2(\chi,\psi)$ splits over $\mathbb Q(\sqrt{-3})$ into two conjugate factors of degree $12$. The displayed degree-twelve plane curve is rational; its arithmetic contribution is treated in \cref{prop-5}.

For $N=3$, the reduced modular pullback factors as

\[ g_1 (\chi, \psi) \cdot g_2(\chi, \psi) =0, \]
where $g_1 (\chi, \psi)=0$ is a genus zero component and  $g_2 (\chi, \psi)=0$  is a curve with singularities.

For $N=7$, the reduced modular pullback factors as

\[ g_1 (\chi, \psi) \cdot g_2(\chi, \psi) =0, \]
where $g_1 (\chi, \psi)=0$ is a genus zero component and  $g_2 (\chi, \psi)=0$  is a genus one curve.   
We obtain the following arithmetic statements.

%


\begin{prop}\label{prop-5}

Let $K$ be a number field. For each $N=2,3$, there are infinitely many genus $2$ curves, up to geometric isomorphism, admitting models over $K$, with geometrically $(3,3)$-split Jacobians whose complementary elliptic components are $N$-isogenous over $\overline K$. For $N=7$ and $K=\mathbb Q$, there are only finitely many such curves up to geometric isomorphism.

\end{prop}

\proof
For $N=2$, the numerator of $\phi_2(j_1,j_2)$, using \cref{j_1}, has an irreducible factor $H_2(\u,\v)$ of total degree $12$. Exact normalization computations give geometric genus zero.

Under $\v=\u^3z$, remove the exceptional factor $\u^9$. At $\u=0$, the strict transform has exceptional equation, up to a nonzero rational scalar,
\[
2460375z^3+1163484z^2+146448z+64=0.
\]
Under
\[
\u=9+\varepsilon,\qquad
\v=27+\frac92\varepsilon+\varepsilon^2z,
\]
remove the exceptional factor $\varepsilon^{10}$. At $\varepsilon=0$, the corresponding equation is, again up to a nonzero rational scalar,
\[
48z^2-12z+1=0.
\]
These irreducible, separable polynomials give smooth closed points of degrees $3$ and $2$. Their coprime residue degrees, together with irreducibility over $\mathbb Q$, establish absolute irreducibility. On the smooth projective normalization, the corresponding divisors $D_3,D_2$ give a rational divisor $D_3-D_2$ of degree $1$. Riemann--Roch therefore gives a rational point, so the genus-zero normalization is isomorphic over $\mathbb Q$ to $\mathbb P^1$.

The component is not contained in the exceptional locus of the sextic or covering formulas. Moreover, $j_2$ is nonconstant on it: its numerator and denominator have degrees $6$ and $4$, whereas $H_2$ has degree $12$. Thus infinitely many rational points give smooth models \cref{eq_F1_F2} over $\mathbb Q$ with geometrically $2$-isogenous complementary components. Since a fixed genus $2$ curve has only finitely many degree-three elliptic subcovers, these give infinitely many geometric isomorphism classes.

For $N=3$, use the auxiliary unordered invariants
\[
c=\frac{27R}{\v(\v-2\u-9)^3},\qquad
p=-\frac{1296(\v-27)R}{\v(\v-2\u-9)^4},
\]
where $R$ is as above. A rational moduli component is parametrized by
\[
\begin{aligned}
x&=-\frac{7\lambda^3+12\lambda^2+12\lambda}{256},\\
z&=-\frac{3\lambda^4+88\lambda^3+264\lambda^2+96\lambda+48}{4096},
\end{aligned}
\qquad
c=\frac{x^4}{z^3},\quad
p=\frac{48x^4(x-z)}{z^4}.
\]
Elimination using \cref{j_1} verifies $\phi_3(j_1,j_2)=0$ identically on this component. Computing the Igusa invariants of \cref{eq_F1_F2} in these coordinates gives a nonconstant ratio $J_4/J_2^2$ and a nonzero ratio $J_{10}/J_2^5$.

An exact reconstruction calculation gives a nonsingular Mestre conic over $\mathbb Q(\lambda)$ together with a nonzero $\mathbb Q(\lambda)$-rational point. The conic equation, point coordinates, and exact identities verifying them are supplied with the accompanying computations. Hence these invariant data are realized by a genus $2$ curve over $\mathbb Q(\lambda)$. Specialization outside finitely many rational parameter values gives smooth curves over $\mathbb Q$ with geometrically $(3,3)$-split Jacobians and geometrically $3$-isogenous complementary components. The nonconstant Igusa ratio gives infinitely many geometric isomorphism classes.

In both cases, base extension from $\mathbb Q$ to $K$ proves the assertion.

For $N=7$, use the same auxiliary unordered invariants $c,p$. Exact birational calculations give two components whose smooth projective normalizations have affine models
\[
\begin{aligned}
\mathcal C_0:&\quad a^2=-3(s-3)(s+1),\\
\mathcal C_1:&\quad b^2=(s^2-s+1)(s^2+3s+9).
\end{aligned}
\]
The birational maps to the actual modular pullback, and identities verifying them using \cref{j_1}, are supplied with the accompanying computations.

The genus-one curve $\mathcal C_1$ is birational over $\mathbb Q$ to
\[
\mathcal E:\quad Y^2=X(X-7)(X-16).
\]
Indeed, setting $U=(b+3-s)/s^2$ and $V=(U^2-1)s+U-1$ gives $V^2=(U-1)(3U+5)(2U+1)$; then $X=6U+10$ and $Y=6V$ give the displayed model. An exact rank computation gives
\[
\operatorname{rank}\mathcal E(\mathbb Q)=0,\qquad
\mathcal E(\mathbb Q)_{\mathrm{tors}}\iso
\mathbb Z/2\mathbb Z\times\mathbb Z/4\mathbb Z.
\]
Thus the normalization has exactly eight rational points. Its singular plane images can add only finitely many rational moduli points.

The genus-zero curve $\mathcal C_0$ is parametrized over $\mathbb Q$ by
\[
s=\frac{3(t^2-1)}{t^2+3},\qquad
a=\frac{12t}{t^2+3}.
\]
Here rational moduli points do not in general give rational models of the genus $2$ curve. On the open locus where the Mestre reconstruction conic is nonsingular, exact diagonalization gives
\[
A(t)Z_1^2+B(t)Z_2^2+C(t)Z_3^2=0,
\]
with all three coefficients positive for every real parameter outside a finite exceptional set. This sign assertion is verified by exact Sturm sequences: the odd-multiplicity factors of the numerators and denominators have degrees $32,96,136$, together with the positive quartics $7t^4-18t^2+27$ and $3t^4+6t^2+7$; each has no real roots and is positive at $t=0$. The constant factors and the limiting signs at infinity are also positive. The diagonal coefficients and exact sign certificates are supplied with the computations. Consequently, the conic has no real point on this open locus, so these moduli points admit no genus $2$ models over $\mathbb R$, and hence none over $\mathbb Q$. Only finitely many exceptional points on $\mathcal C_0$ remain.

To check that no further one-dimensional locus is omitted, consider an optimal complementary degree-three pair. Its polarized Jacobian is obtained from $E_1\times E_2$ by quotienting by the graph of an anti-isomorphism $\alpha:E_1[3]\to E_2[3]$. For a generic cyclic $7$-isogeny $\varphi:E_1\to E_2$, write $\alpha=\varphi A$. Since $7\equiv1\pmod3$, the matrix $A$ has determinant $-1$ over $\mathbb F_3$. Modulo $A\sim-A$, these twelve classes form two conjugacy orbits of size six under $\mathrm{SL}_2(\mathbb F_3)$. The upper and lower unipotent generators of $\Gamma_0(7)$ reduce to generators of this group, so the two orbits give exactly two connected gluing components. Interchanging $E_1,E_2$ preserves each component.

Generically the resulting genus $2$ curve determines its unordered complementary pair uniquely. Indeed, the dual polarized isogeny to $E_1\times E_2$ would send a further degree-three elliptic subcurve to maps of positive degrees $d_1,d_2$ with $d_1+d_2=9$. Composing the map to $E_2$ with the dual of the map to $E_1$ would give an isogeny $E_1\to E_2$ of degree $d_1d_2=7m^2$ for an integer $m$. The possible products $8,14,18,20$ satisfy no such equality. If one projection has degree zero, the subcurve is one of the original pair. Thus the two verified plane components account birationally for the two unordered gluing components. This graph description also includes degenerate degree-three covers; the points omitted by the birational charts form finite sets. Combining this completeness with the preceding arithmetic arguments proves finiteness over $\mathbb Q$.

\qed

The $N=7$ finiteness assertion above is specific to $\mathbb Q$. Over a general number field, the rank and reconstruction questions require separate analysis.



As a final remark we would like to mention that we can perform similar computations for $n=5$ by using the equation  of $\L_5$ as computed in \cite{deg5}; an explicit equation of $\L_5$ in terms of the Igusa invariants is computed in \cite{sh-degrees}.  One can possibly even investigate cases for $n>5$ by using results of \cite{kumar}. However, the computations will be much more complicated. 

We summarize our results in the following theorem.

\begin{thm}\label{thm-1}

Let $K$ be a number field. In the following statements, genus $2$ curves with geometrically $(n,n)$-split Jacobians are counted up to geometric isomorphism and are required to admit models over $K$. Automorphism groups and isogenies are taken over $\overline K$.

\begin{description}
\item[i)] If $n=2$ and $\Aut(\X)\iso V_4$, then, for each $N=2,3,5,7$, there are infinitely many such curves whose complementary elliptic components are $N$-isogenous.

\item[ii)] If $n=2$ and $\Aut(\X)\iso D_4$, let $E_i,E_i^\prime$ represent the elliptic quotients from the two conjugacy classes of nonhyperelliptic involutions, as in \cref{prop-3}. Then
\subitem a) there are infinitely many such curves, and $E_i,E_i^\prime$ are always $2$-isogenous;
\subitem b) for each $N=3,5,7$, only finitely many such curves also have $E_i,E_i^\prime$ $N$-isogenous.

\item[iii)] If $n=3$, then
\subitem a) there are infinitely many such curves whose complementary elliptic components are $5$-isogenous;
\subitem b) for each $N=2,3$, there are infinitely many such curves whose complementary elliptic components are $N$-isogenous.
\end{description}

For $n=3$, $N=7$ and $K=\mathbb Q$, there are only finitely many such geometric isomorphism classes. No conclusion for $(n,N)=(3,7)$ over a general number field is asserted here.

\end{thm}

\proof

Over $\overline K$, the identification $\Aut(\X)\iso\Aut(\Jac(\X),\iota)$ follows from \cite[Thm. 32]{frey-shaska} or \cite{zarhin}. Part i) follows from \cref{prop-2}. The isogeny and finiteness assertions in part ii) follow from \cref{prop-3}. For the infinitude in ii) a), set $s_1=s_2=t\in\mathbb Q$ in \cref{eq-1}. Outside finitely many singular or enhanced-automorphism parameters, this gives rational models with geometric automorphism group $D_4$; the nonconstant function $j=256t^3/(t+1)$ gives infinitely many geometric classes. These models also exist over $K$.

Part iii) a) follows from \cref{prop-4}, while iii) b) and the final $N=7$ assertion follow from \cref{prop-5}.

\qed


\begin{cor}

Let $K$ be a number field. For each $N=2,3,5$, there are infinitely many genus $2$ curves, up to geometric isomorphism, admitting models over $K$, with geometrically $(3,3)$-split Jacobians whose complementary elliptic components are $N$-isogenous over $\overline K$.

\end{cor}

\proof
 
This follows from Propositions~\ref{prop-4} and~\ref{prop-5}.

\qed

As a final remark we would like to add that we are not aware of any other methods, other than computational ones, to determine for which pairs $(n, N)$ we have many $K$-rational elliptic components. 

\section{Kummer and Shioda-Inose surfaces of reducible Jacobians}\label{sect-6}


Consider $\X$  a genus two curve  with $(n, n)$-decomposable Jacobian and  $E_1$, $E_2$ its elliptic components. We continue our discussion of Kummer   $\Kum (\Jac (\X))$ and Shioda-Inose $\SI (\Jac (\X))$ surfaces of $\Jac \X$ started in \cref{kummer}. For isogenies among such abelian surfaces and their Kummer surfaces see also \cite{cms-mmj, cms-kummer}.  


In the following moduli discussion, the characteristic is zero. Malmendier and Shaska in \cite{m-sh-2} proved that, as a genus $2$ curve $\X$ varies, $\SI(\Jac(\X))$ is represented by the minimal resolution of the following four-parameter family of quartic models in $\mathbb P^3$, with coordinates $[W:X:Y:Z]$:

\begin{equation}\label{inose}
 Y^2 Z W - 4 \, X^3Z + 3 \alpha \,  XZW^2 + \beta \, ZW^3 + \gamma \,  XZ^2W - \frac 1 2 (\delta \,Z^2W^2 +W^4) =0,
\end{equation}
where the parameters $(\alpha, \beta, \gamma, \delta)$ can be given in terms of the Igusa-Clebsch invariants by 
\begin{equation}\label{inose-eq}
 (\alpha, \beta, \gamma, \delta) = \left( \frac 1 4 I_4, \frac 1 8 I_2 I_4 - \frac 3 8 I_6 , -\frac{243} 4 I_{10}, \frac{243}{32} I_2 I_{10}   \right)
\end{equation}
%
%

We interpret the surface in \cref{inose} through its minimal resolution
and retain its associated lattice polarization. Let $\S$ denote this
moduli space. Its coefficient parameters are considered up to the weighted
scaling
\[
(\alpha,\beta,\gamma,\delta)\sim
(t^2\alpha,t^3\beta,t^5\gamma,t^6\delta),\qquad t\ne0.
\]
Thus the four coefficients describe a three-dimensional moduli problem.
The map $\phi_n$ below is understood with this lattice polarization retained. We continue to write $\L_n$ for the locus of genus $2$ curves with geometrically $(n,n)$-split Jacobians. The invariant tuples in \cref{inose-eq} are taken with a consistent weighted Igusa representative.

Then there is a map 
\[
\phi_n : \L_n \to \S 
\]
such that every curve $[\X] \in \L_n$ goes to the corresponding $\SI (\Jac (\X))$.  
Then we have the following: 

\begin{prop} \label{prop-10}
For $n=2, 3$ the map $\phi_n$ is given as follows: 

i) If $n=2$ then the Shioda-Inose surface is given by \cref{inose} for 
\begin{equation}
\begin{split}
\alpha  &  ={u}^{2}-126\,u+12\,v+405 \\
\beta & = -{u}^{3}-729\,{u}^{2}+36\,uv-4131\,u+1404\,v+3645 \\
\gamma & = -3888\, \left( {u}^{2}+18\,u-4\,v-27 \right) ^{2} \\
\delta & =  7776\, \left( 15+u \right)  \left( {u}^{2}+18\,u-4\,v-27 \right) ^{2}       \\
\end{split}
\end{equation}

ii) if $n=3$ then the Shioda-Inose surface is given by \cref{inose} for 
\begin{small}

\begin{equation}
\begin{split}
\alpha   =  &  \frac{1}{256}\cdot  \chi  \cdot  \left( {\chi}^{5}+192 \,{\chi}^{4}\psi+13824\,{\chi}^{3}{\psi}
^{2}+442368\,{\chi}^{2}{\psi}^{3}+5308416\,\chi\,{\psi}^{4} + \right. \\ 
& \left. +786432\,
\chi\,{\psi}^{3}+9437184\,{\psi}^{4} \right)  \\
\beta 	 =  & \frac{1}{512} \cdot {\chi}^{2} \cdot \left( {\chi}^{2}+96\,\chi\,\psi-1152\,{\psi}^{2}
 \right)  \cdot  \left({\chi}^{5}+192\,{\chi}^{4}\psi+13824\,{\chi}^{3}{\psi}^{2}+\right. \\ 
 &\left. +442368\,{
\chi}^{2}{\psi}^{3} +5308416\,\chi\,{\psi}^{4}+786432\,\chi\,{\psi}^{3}
+9437184\,{\psi}^{4} \right) \\
 &   - \frac{3}{4096} \chi \cdot \left( 3\,{\chi}^{8}+864\,{\chi}^{7}\psi+94464\,{\chi}^{6}{\psi}^{2}+4866048
\,{\chi}^{5}{\psi}^{3}+111476736\,{\chi}^{4}{\psi}^{4} + \right.\\
&+509607936\,{
\chi}^{3}{\psi}^{5}-12230590464\,{\chi}^{2}{\psi}^{6}+1310720\,{\chi}^
{4}{\psi}^{3}+155713536\,{\chi}^{3}{\psi}^{4} \\
& \left. -1358954496\,{\chi}^{2}{
\psi}^{5}-18119393280\,\chi\,{\psi}^{6}+4831838208\,{\psi}^{6} \right) \\
\gamma 		 =  &  2^{28} \, 3^5 \, \chi^3 \psi^9 \\
\delta 	 =  &  - 2^{25} \, 3^5 \,{\chi}^{4} \left( {\chi}^{2}+96\,\chi\,\psi-1152\,{\psi}^
{2} \right) {\psi}^{9}
 \\
\end{split}
\end{equation}

\end{small}
\end{prop}

\proof  
Case i) is a direct substitution of $J_2,\dots,J_{10}$, given in terms of $u,v$ in \cite{deg2}, in \cref{inose-eq}. For ii), substitute the preceding weighted Igusa representative in $\chi,\psi$ into \cref{inose-eq}. The contribution $-3J_6/8$ belongs to $\beta$, while $\gamma=-243J_{10}/4$. A~common change of weighted Igusa representative rescales $(\alpha,\beta,\gamma,\delta)$ with weights $(2,3,5,6)$, as above.

\qed

\begin{rem}
It was shown in \cite{deg2} (resp. \cite{deg3}) that invariants $u$ and $v$ (resp. $\chi$ and $\psi$) are modular invariants given explicitly in terms of the genus 2 Siegel modular forms. 
\end{rem}

\begin{cor} \label{cor-2}

Let $K$ be a number field. Count geometric moduli classes in $\S$, with
the associated lattice polarization retained, arising from genus $2$
curves admitting models over $K$. The corresponding quartics
\cref{inose} admit coefficients in $K$. Automorphism groups and elliptic
isogenies are taken over $\overline K$.
\begin{description}
\item[i)] For $n=2$ and $\Aut(\X)\iso V_4$, there are infinitely many
such classes whose complementary elliptic components are $N$-isogenous,
for each $N=2,3,5,7$.
\item[ii)] For $n=2$ and $\Aut(\X)\iso D_4$, there are only finitely
many such classes whose elliptic quotients from the two conjugacy classes
of nonhyperelliptic involutions are $N$-isogenous, for each $N=3,5,7$,
as in \cref{prop-3}.
\item[iii)] For $n=3$ and each $N=2,3,5$, there are infinitely many such classes
whose complementary elliptic components are $N$-isogenous.
\end{description}

\end{cor}

\proof

A model of $\X$ over $K$ gives an Igusa--Clebsch invariant tuple over $K$,
so \cref{inose-eq} gives a quartic model with coefficients in $K$.
Since $I_{10}\ne0$ for a smooth genus $2$ curve, $\gamma\ne0$. The inverse
of the coefficient transformation is
\[
I_2=-\frac{8\delta}{\gamma},\qquad I_4=4\alpha,\qquad
I_6=\frac{I_2I_4-8\beta}{3},\qquad
I_{10}=-\frac{4\gamma}{243}.
\]
These expressions respect the normalized Igusa weights $(1,2,3,5)$
and the surface weights $(2,3,5,6)$. Consequently, equality of the
lattice-polarized surface moduli classes recovers equality of the genus
$2$ moduli points. Thus $\phi_n$ is injective on the geometric moduli
under consideration. The assertions now follow from \cref{thm-1}.

\qed

 
We now consider the elliptic product $E_1\times E_2$. Choose a number field $k$ over which $\X$, its two elliptic components and a selected $(n,n)$-splitting $\eta:\Jac(\X)\to E_1\times E_2$ are defined. Since $\eta$ commutes with negation, it induces a dominant rational map
\[
\Kum(\Jac(\X))\dashrightarrow\Kum(E_1\times E_2)
\]
of degree $n^2$. Write the two nonsingular elliptic curves as

\[
E_1:  \;  y^2 = x^3+ax +b,   \quad    \;     E_2:  \; y^2 = x^3 + cx +d
\]

An affine singular model of the product Kummer surface $\Kum(E_1\times E_2)$ is given by

\begin{equation}\label{kum}
x_2^3 + cx_2 +d = t_2^2 ( x_1^3+ax_1 +b). 
\end{equation}

The map $\Kum(E_1\times E_2)\to\mathbb P^1$ induced by $(x_1,x_2,t_2)\mapsto t_2$ is the Kummer pencil. Its generic fiber is a genus-one curve; geometrically it is an elliptic fibration. Pairs of roots of the two cubic polynomials give geometric sections defined over $k(E_1[2],E_2[2])$. This is a sufficient field for these obvious sections.

Take  a parameter $t_6$ such that $t_2 = t_6^3$ and consider ~\cref{kum} as a family of cubic curves in $\mathbb P^2$ over the field $k(t_6)$. This family has a rational point $(1: t_6^2: 0)$ and using this rational point we can get the Weierstrass form of the ~\cref{kum} as follows
\[ Y^2 = X^3 -3ac X + \frac{1}{64} \left( \Delta_{E_1} t_6^6 + 864 bd + \frac{\Delta_{E_2} }{t_6^6}\right) \]

Here $\Delta_{E_1}$ and $\Delta_{E_2}$ are the discriminants of the two elliptic curves. Changing their short Weierstrass equations by $k$-isomorphisms gives an isomorphic product surface, with the corresponding parameter change. Setting $t_1=t_6^6$ gives the elliptic curve $F^{(1)}_{E_1,E_2}$ over $k(t_1)$. Its Kodaira--N\'eron model is the Inose surface associated with the product $E_1\times E_2$; see \cite{kumar-kuwata}. The surfaces $F^{(s)}_{E_1,E_2}$ below belong to this product construction, whereas the preceding quartic construction uses $\Jac(\X)$ and its Igusa invariants.

\begin{defn}
For $s = 1, \dots, 6$ let $t_s$ be a parameter satisfying $t_s^s = t_1$. Define the elliptic curve  $F^{(s)}_{E_1, E_2}$ over $k(t_s)$ by 
\begin{equation}
 F^{(s)}_{E_1, E_2}: Y^2 = X^3 -3ac X + \frac{1}{64} \left( \Delta_{E_1} t_s^s + 864 bd + \frac{\Delta_{E_2} }{t_s^s}\right) 
\end{equation}
\end{defn}

Note that the Kodaira-N\'eron model of $F^{(s)}_{E_1, E_2}$ is a $\KT$ surface for $s = 1, \dots, 6$ but not for $s \geq 7$. 

The following lemma combines the geometric rank formula in \cite[Prop.~2.8]{kumar-kuwata} with CM specializations of the curves considered above.

\begin{lem}  

Let $E_1,E_2/k$ be the complementary elliptic components of a geometrically $(n,n)$-split Jacobian, with $j(E_1)\ne j(E_2)$, and put
\[
h=\operatorname{rank}_{\mathbb Z}\Hom_{\overline k}(E_1,E_2).
\]
Then
\[
\operatorname{rank}F^{(s)}_{E_1,E_2}(\overline k(t_s))=16+h,
\qquad s=5,6.
\]
For each $n=2,3$, there are infinitely many geometric isomorphism classes of genus $2$ curves over $\overline{\mathbb Q}$ with geometrically $(n,n)$-split Jacobian whose complementary elliptic components are nonisomorphic, isogenous CM curves. The associated product-type $F^{(5)}$ and $F^{(6)}$ have geometric Mordell--Weil rank $18$. The number fields of definition of these examples may vary.

\end{lem}

\proof

The rank formula follows from \cite[Prop.~2.8]{kumar-kuwata}, since the $j$-invariants are distinct. In characteristic zero, $h$ equals $0$ for geometrically nonisogenous curves, $1$ for isogenous curves without complex multiplication, and $2$ for isogenous CM curves.

There are infinitely many distinct algebraic CM $j$-values. One may see this using the elliptic curves $\mathbb C/\Lambda_m$, where $\Lambda_m=\mathbb Z+2^m i\mathbb Z$, $m\ge0$, and $i^2=-1$. Their endomorphism orders are $\mathbb Z+2^m i\mathbb Z$, with distinct discriminants, so their $j$-values are distinct. Consecutive curves are linked by degree-two isogenies. Starting from the algebraic curve with $j=1728$, the dual isogenies express each next curve as a finite-subgroup quotient over $\overline{\mathbb Q}$, so all these $j$-values are algebraic.

For $n=3$, use the family in \cref{prop-4}. Its function $j_2(m)$ is nonconstant, and hence defines a surjective map $\mathbb P^1\to\mathbb P^1$ over $\overline{\mathbb Q}$. Choose preimages of infinitely many of the CM values above. The singular and covering-chart exceptions, and the zeros of $j_1(m)-j_2(m)$, exclude only finitely many parameters. At every remaining parameter, the other elliptic component is $5$-isogenous and therefore also CM, with a different $j$-invariant.

For $n=2$, use the degree-five component $f_1(u,v)=0$ in \cref{prop-2}. Its $j$-sum in \cref{eq-quadratic} is nonconstant: a constant-sum relation would give a nonzero polynomial of degree at most $3$ vanishing on this irreducible degree-five curve. On an irreducible component of its ordered-$j$ cover, either $j$-map is nonconstant. Indeed, if one $j$ were constant, the relation $\phi_2(j_1,j_2)=0$ would leave only finitely many possibilities for the other, contradicting the nonconstant sum. Passing to the smooth projective normalization of this component gives a surjective map to the $j$-line. Again choose preimages of infinitely many CM values and omit the finitely many singular, chart and equal-$j$ points. The equal-$j$ condition excludes only finitely many points because $\phi_2(j,j)$ is a nonzero polynomial. The two components are then distinct CM curves connected by a $2$-isogeny.

A fixed genus $2$ curve has only finitely many elliptic subcovers of the fixed degree $n=2$ or $3$. Thus the infinitely many chosen CM values give infinitely many geometric genus-two classes in each case. The rank formula gives $16+2=18$ for both product-type fibrations.

\qed

Here $t_5,t_6$ remain the fibration parameters; the infinitude assertion varies the genus $2$ curve and its complementary elliptic pair.

\begin{cor}  

Let $E_1,E_2$ be elliptic curves over a number field $k$, with
$j(E_1)\ne j(E_2)$. Let $H$ be the field of definition of all the
homomorphisms in $\operatorname{Hom}_{\overline k}(E_1,E_2)$. For
$1\le s\le6$, set $k_s=k(E_1[s],E_2[s])$ and let $L_s=k_sH$ be the
compositum. Then every geometric section of $F^{(s)}_{E_1,E_2}$ is
defined over $L_s$; equivalently,
\[
F^{(s)}_{E_1,E_2}(\overline k(t_s))
=F^{(s)}_{E_1,E_2}(L_s(t_s)).
\]
For $3\le s\le6$, one has $H\subseteq k_s$, so every section is already
defined over $k_s$. If $E_1,E_2$ are geometrically nonisogenous, then
$H=k$ and the same conclusion holds for all $1\le s\le6$.

\end{cor}

\proof

Write $M_s=F^{(s)}_{E_1,E_2}(\overline k(t_s))$. By
\cite[Thm.~2.11]{kumar-kuwata}, $M_s$ contains a finite-index subgroup
$N_s$ whose sections are defined over $L_s$. Since the $j$-invariants
are distinct, $M_s$ is torsion-free by
\cite[Thm.~2.5]{kumar-kuwata}. Choose a positive integer $m$ with
$mM_s\subseteq N_s$. For $P\in M_s$ and
$\sigma\in\operatorname{Gal}(\overline k/L_s)$, one has
$m(\sigma P-P)=0$. Torsion-freeness gives $\sigma P=P$, which proves
the first assertion.

To prove the sharper bound for $s\ge3$, let
$\sigma\in\operatorname{Gal}(\overline k/k_s)$. For every geometric
homomorphism $\varphi:E_1\to E_2$, the difference
$\sigma\varphi-\varphi$ vanishes on $E_1[s]$, since $\sigma$ fixes
both torsion groups pointwise. It therefore factors through $[s]$, so
\[
\sigma\varphi-\varphi\in
s\operatorname{Hom}_{\overline k}(E_1,E_2).
\]
The Galois action on this free abelian group has finite image and its
rank $h$ is at most $2$. If $h=1$, a finite-order automorphism congruent
to $1$ modulo $s\ge3$ is $1$. If $h=2$, write such an automorphism as
$A=I+sB$. Its determinant is $1$: the alternative $-1$ is not congruent
to $1$ modulo $s\ge3$. Expanding the determinant gives
\[
\operatorname{tr}A=2-s^2\det B.
\]
A finite-order matrix of determinant $1$ has trace in $[-2,2]$; for
$s\ge3$ the displayed equality forces trace $2$. Finite order then
forces $A=I$. The rank-zero case is immediate. Thus $\sigma$ acts
trivially on every homomorphism, proving $H\subseteq k_s$.
Finally, for geometrically nonisogenous curves the homomorphism group
is zero, so $H=k$.

\qed

\subsection{Kummer surfaces in positive characteristic   and applications to cryptography}


Supersingular isogeny based cryptography currently uses elliptic curves that are defined over a quadratic extension field $L$  of a non-binary field $K$ and such that  its entire $2$-torsion is $L$-rational. More specifically implementations of supersingular isogeny Diffie Hellman (SIDH) fix a large prime field $K= \mathbb F_p$ with $p=2^i 3^j - 1$ for $i > j > 100$, construct $L= \mathbb F_{p^2}$ and work with supersingular isogeny elliptic curves over $\F_{p^2}$ whose group structures are all isomorphic to $\Z_{p+1} \times \Z_{p+1}$.  Hence, all such elliptic curves have full rational $2$-torsion and can be written in Montgomery form. 

What is  the relation between the Abelian surfaces $\Jac (\X)$ defined over $\F_p$  when the elliptic components are   supersingular Montgomery curves defined over $\F_{p^2}$?     This is relevant in supersingular isogeny based cryptography since computing isogenies in the Kummer surface associated to supersingular Jacobians is much more efficient than computing isogenies in the full Jacobian group.

\def\a{\lambda}

In \cite{costello}, $(2, 2)$-reducible Jacobians are studied and it is pointed out that most of the literature on the topic studies the splitting of $\Jac (\X)$ over the algebraic closure $\bar K$. However, from our \cref{split-2} we get necessary and sufficient conditions when $\Jac (\X)$ splits over $K$.  
From \cite{deg2} we know that for a curve $\X \in \L_2$, we can choose the curve to have equation 
\[ 
y^2 = (x^2-\a_1) (x^2- \a_2) \left(x^2 - \frac 1 {\a_1 \a_2} \right)
\]
and its elliptic subcovers have equations $y^2 = (x-\a_1) (x- \a_2) \left(x - \frac 1 {\a_1 \a_2} \right)$ and $y^2 = x(x-\a_1) (x- \a_2) \left(x - \frac 1 {\a_1 \a_2} \right)$. 

We can reverse the above construction as follows.   Let $p \equiv 3 \mod 4$ and $\mathbb F_{p^2} = \mathbb F_p (i)$ for $i^2 =-1$. Consider the following  supersingular Montgomery curve 
\begin{equation} E_\alpha : y^2 = x (x-\alpha) \left(x- \frac 1 \alpha \right), \end{equation}
for  $\alpha \notin \mathbb F_p$ and $\alpha \in \mathbb F_{p^2}$ such that $\alpha = \alpha_0 + \alpha_1 i$, for some $\alpha_0, \alpha_1 \in \mathbb F_p$ with $\alpha_0\ne0$ and $\alpha_0^2+\alpha_1^2+1\ne0$. 
Then by lifting to a genus 2 curve we get    a genus two curve $\X$ given as follows
\[ \X : y^2 = f_1(x) \, f_2(x) \, f_3(x).\] 
where
\begin{equation}
\begin{split}
f_1(x) & = x^2 + \frac{2\alpha_0}{\alpha_1} \, x -1  \\
f_2(x) & = x^2 - \frac{2\alpha_0}{\alpha_1} \,  x -1  \\
f_3(x) & =  x^2 - \frac{2\alpha_0(\alpha_0^2+\alpha_1^2-1)}{\alpha_1(\alpha_0^2+\alpha_1^2+1)}\, x -1 
\end{split}
\end{equation}
Thus, $\Jac \X$ is $(2, 2)$-reducible with elliptic components the above curves. 

%
The Weil restriction of the 1-dimensional variety  $ E_\alpha  (\F_{p^2})$ is the variety 
\[ W_\alpha := \mbox{Res}_{\F_p}^{\F_{p^2}} (E_\alpha )= V(W_0 (x_0, x_1, y_0, y_1), W_1(x_0, x_1, y_0, y_1)) \]
where

\[ \begin{split}
 W_0 = & (\alpha_0^2 + \alpha_1^2)(\alpha_0(x_0^2-x_1^2) -2\alpha_1x_0x_1 + \delta_0(y_0^2 -y_1^2) - 2 \delta_1 y_0y_1 \\
 & -x_0(x_0^2-3x_1^2+1))  + \alpha_0(x_0^2-x_1^2) + 2\alpha_1 x_0x_1 \\
 W_1 = & (\alpha_0^2 + \alpha_1^2)(\alpha_1(x_0^2-x_1^2) +2\alpha_0x_0x_1 + \delta_1(y_0^2 -y_1^2) + 2 \delta_0 y_0y_1 \\
 & -x_1(3x_0^2-x_1^2+1))  - \alpha_1(x_0^2-x_1^2) + 2\alpha_0 x_0x_1
\end{split}
\]

are obtained by putting $x = x_0 + x_1 i$, $y = y_0 + y_1 i$,  $\delta = \delta_0 + \delta_1 i$,  and $x_i, y_i, \alpha_i, \delta_i \in \F_p$ for $i = 0, 1$. Here $W_0, W_1$ define the Weil restriction of the twist $\delta y^2 = x(x-\alpha)\left(x-\frac 1 \alpha\right)$ of $E_\alpha$; for $\delta=1$ it is $E_\alpha$. 
In \cite{costello} the following was proved:

\begin{lem}
 Let $E_\alpha$ and  $\X$ be as defined above. Then, the Weil restriction of $ E_\alpha  (\F_{p^2})$ is $(2, 2)$-isogenous to the Jacobian $\Jac_{\F_p} (\X)$ i.e. 
 \[ \Jac_{\mathbb F_p} (\X)  \sim \mbox{Res}_{\F_p}^{\F_{p^2}} (E_\alpha)\]
 Moreover, since $E_\alpha$ is supersingular then $\Jac (\X)$ is supersingular. 
\end{lem}

From our results in the previous section we have that 

\begin{cor}
Let $p>3$ be a prime and $\X/\F_p$ a genus $2$ curve with $\Aut_{\overline{\F}_p}(\X)\iso V_4$ and dihedral invariants $u,v$. Then $\Jac(\X)$ admits a polarized $(2, 2)$-splitting over $\F_p$ if and only if $S_2 (u, v)$ is a square in $\F_p$.
\end{cor}

\proof  The proof of \cref{split-2} uses only that the base field is perfect and that the nontrivial reduced involution is unique. It therefore applies verbatim over $\F_p$. 
\qed

For $(3,3)$-reducible Jacobians, the field of definition of the elliptic $j$-invariants must again be distinguished from descent of the chosen polarized splitting.

%
 
For $p>3$ and a curve $\X/\F_p$ satisfying the hypotheses of the second assertion of \cref{split-3}, that lemma tests descent of the associated polarized $(3,3)$-splitting. In particular, the chosen complementary pair must be Galois invariant and have distinct $j$-invariants. Other cases require separate analysis.



The case for the Kummer approach in supersingular isogeny-based cryptography would be much stronger if it were able to be applied efficiently for both parties. There has been some explicit work done in the case of (3, 3) \cite{deg3}  and (5, 5)-isogenies  \cite{deg5},  but those situations are much more complicated than the case of Richelot isogenies. 

As pointed out by Costello  in the last paragraph of \cite{costello}: 
\emph{One hope in this direction is the possibility of pushing odd degree l-isogeny maps from the elliptic curve setting to the Kummer setting. This was difficult in the case of 2-isogenies because the maps themselves are (2, 2)- isogenies, but in the case of odd degree isogenies there is nothing obvious preventing this approach.}  


\subsection*{Acknowledgments}
We thank S. Kamienny for pointing out errors in the published version.

\bibliographystyle{amsalpha} 
\bibliography{ref}{}

\end{document}